\documentclass[reqno]{article}

\usepackage{xcolor}
\usepackage{amssymb}
\usepackage{amsfonts}
\usepackage{amsmath}
\numberwithin{equation}{section}
\usepackage{enumerate}

\newtheorem{theorem}{Theorem}[section]
\newtheorem{proposition}[theorem]{Proposition}
\newtheorem{lemma}[theorem]{Lemma}
\newtheorem{corollary}[theorem]{Corollary}
\newtheorem{definition}[theorem]{Definition}
\newtheorem{example}[theorem]{Example}

\newtheorem{remark}[theorem]{Remark}

\newcommand{\ran}{{\rm ran}}

\newcommand{\dom}{{\rm dom}}
\newcommand{\mul}{{\rm mul}}

\title {General Green's boundary model in a Krein space}
 
\author{Author: Muhamed Borogovac}

\begin{document}

\maketitle 

%\textbf{A concise and informative title}:  Green's boundary relation
%\\
%\textbf{The affiliation and address of the author}: 

BML - Actuarial Department, Canton, MA 02021, USA
%\textbf{The e-mail address, and telephone of the author}: 

muhamed.borogovac@gmail.com
%1-857-991-8779 

The 16-digit ORCID of the author: \textbf{0000-0001-6857-8309}
\\

\textbf{Abstract:} Given Krein and Hilbert spaces $\left( \mathcal{K},[.,.] \right)$ and $\left( \mathcal{H}, \left( .,. \right) \right)$, respectively, the concept of the boundary triple $\Pi =(\mathcal{H}, \Gamma _{0}, \Gamma_{1})$ is generalized through the abstract Green's identity for the isometric relation $\Gamma$ between Krein spaces $\left( \mathcal{K}^{2}, \left[ .,.\right]_{\mathcal{K}^{2}} \right) $ and $\left(\mathcal{H}^{2}, \left[ .,.\right]_{\mathcal{H}^{2}} \right) $ without any conditions on $\dom\, \Gamma$ and $\ran\, \Gamma$. This also means that we do not assume the existence of a closed symmetric linear relation $S$ such that  $\dom\, \Gamma=S^{+}$, which is a standard assumptions in all previous research of boundary triples. The main properties of such a general Green's boundary model are proven. In the process, some useful properties of the isometric relation $V$ between two Krein spaces $X$ and $Y$ are proven. Additionally, surprising properties of the unitary relation $\Gamma : \mathcal{K}^{2} \rightarrow\mathcal{H}^{2}$ and the self-adjoint main transformation $\tilde{A}$ of $\Gamma$ are discovered. Then, two statements about generalized Nevanlinna families are generalized using this Green's boundary model. Furthermore, several previously known boundary triples involving a Hilbert space $\mathcal{K}$ and reduction operator $\Gamma : \mathcal{K}^{2} \rightarrow\mathcal{H}^{2}$, such as AB-generalized, B-generalized, ordinary, isometric, unitary, quasi-boundary, and S-generalized boundary triples, have been extended to a Krein space $\mathcal{K}$ and linear relation $\Gamma$ using the Green's boundary model approach.



\section{Preliminaries and introduction }\label{s2}

\textbf{1.1.} By $\mathbb{N}$, $\mathbb{R}$, $\mathbb{C}$ we denote the sets of positive integers, 
real numbers, and complex numbers, respectively. Let $\left( \mathcal{K},\, \left[ .,. 
\right] \right)$ denote a Krein space. That is a complex vector space on 
which a scalar product, i.e., a Hermitian sesquilinear form $\left[ .,. 
\right]$, is defined such that the following decomposition of $\mathcal{K}$ exists: 
\[
\mathcal{K}=\mathcal{K}_{+}[+]\mathcal{K}_{-},
\]
where $\left( \mathcal{K}_{+},[.,.] \right)$ and $\left( \mathcal{K}_{-},-[.,.] \right)$ are Hilbert spaces, and the sum $[+]$ is direct and orthogonal with respect to the form $[.,.]$. We say that elements $x,y \in \mathcal{K} $ are \textit{orthogonal} if $\left[ x,y\right]=0$, denoted by $x \left[ \perp \right]y$. Every Krein space $\left( \mathcal{K},[.,.] \right)$ is \textit{associated} with the Hilbert space $\left( \mathcal{K},(.,.) \right)$, which is defined as a direct and orthogonal sum of the Hilbert spaces $\left( \mathcal{K}_{+},[.,.] \right)$ and $\left( \mathcal{K}_{-},-[.,.] \right)$. The Topology in the Krein space $\mathcal{K}$ is introduced by means of the associated Hilbert space $\left( \mathcal{K},(.,.) \right)$. The \textit{orthogonal companion} $A^{\left[ \perp \right]}$  of the set $A$ is defined by $A^{\left[ \perp \right]}:=\left\lbrace y\in \mathcal{K}:\left[x,y \right]=0, \forall x \in A\right\rbrace $, and the \textit{isotropic} manifold $M$ of $A$ is defined by  $M:= A \cap A^{\left[ \perp \right]}$. A set $A \subseteq \mathcal{K}$ is called \textit{neutral} if it satisfies $A \subseteq A^{\left[ \perp \right]}$. A set $A \subseteq \mathcal{K}$ is called \textit{hyper-maximal neutral} if $A = A^{\left[ \perp \right]}$. For properties of Krein spaces, one can refer to e.g., \cite[Chapter V]{Bog}. We will frequently deal with the well known Krein spaces $\mathcal{K}^{2}=\left( \mathcal{K}^{2},J_{\mathcal{K}}\right) $ and $\mathcal{H}^{2}=\left( \mathcal{H}^{2},J_{\mathcal{H}}\right)$, but for the convenience of the reader, the definitions are repeated in Section \ref{s4}, where we need them. 

If the scalar product $\left[ .,. \right]$ has $\kappa \in \mathbb{N}$ negative squares, then we call it a \textit{Pontryagin space of negative index} $\kappa $. If $\kappa =0$, then it is a Hilbert space. More information about Pontryagin space can be found e.g. in \cite {IKL}. 
\\

\textbf{1.2.} The following definitions of a linear relation and basic concepts related to it can be found in \cite {A,S,DS}. In the sequel, $X$, $Y$, $W$ are Krein spaces which include Pontryagin and Hilbert spaces. 

A \textit{linear relation} $T: X\rightarrow Y$ is a linear manifold $T \subseteq X \times Y $. If $ X=Y $, then $T$ is said to be a \textit{linear relation in} $X$. A linear relation $T$ is closed if it is a (closed) subspace with respect to the product topology of $X \times Y $. As usual, for a linear relation $ T: X \rightarrow Y $, or $ T \subseteq X \times Y $, the symbols $ \dom \, T $, $ \ran\, T $, and $ \ker T $ mean domain, range, and kernel, respectively. In addition, we will use the following concepts and notation for two linear relations, $T$ and $S$ from $X$ into $Y$, and a linear relation $U$ from $Y$ into $W$: 
\[
\mul\, T:=\left\{ g\in Y : \left\{ 0,g \right\}\in T \right\};
\]
\[
T\left( f \right):=\left\{ g\in Y : \left\{ f,g \right\}\in T \right\}, f\in D\left( T \right);
\]
\[
T^{-1}:=\left\{ \left\{ g,f \right\}\in Y \times X : \left\{ f,g \right\}\in T \right\};
\]
\[
zT:=\left\{ \left\{ f,zg \right\}\in X \times Y : \left\{ f,g \right\}\in T \right\}, z\in \mathbb{C};
\]
\[
S+T:=\left\{ \left\{ f,g+k \right\} : \left\{ f,g \right\}\in S,\left\{ f,k \right\}\in T \right\};
\]
\[
S\hat{+}T:=\left\{ \left\{ f+h,g+k \right\} : \, \left\{ f,g \right\}\in S,\left\{ h,k \right\}\in T \right\};
\]
\[
S\dot{+}T:=\left\{ \{f+h,g+k\} : \left\{ f,g \right\}\in S,\left\{ h,k \right\}\in T, S\cap T=\lbrace 0\rbrace \right\};
\]
\[
UT:=\left\{ \left\{ f,k \right\}\in X \times W \, : \, \left\{ f,g \right\}\in T,\left\{ g,k \right\}\in U\, for\, some\, g\in Y \right\};
\]
\[
T^{\ast}:=\left\{ \left\{ k,h \right\}\in Y \times X : \left[ f,h \right]=\left[ g,k \right]\, for\, all\, \left\{ f,g \right\}\in T \right\};
\]
\[
T_{\infty }:=\left\{ \left\{ 0,g \right\}\in T \right\}=\lbrace 0 \rbrace \times \mul\, {T}.
\]
If $\mul\, {T}=\{ 0\}$, we say that T is \textit{single-valued} linear relation, i.e., an \textit{operator}. If additionally $\ker\, {T}=\lbrace 0\rbrace $, then we call it \textit{one-to-one} operator. The sets of closed linear relations, closed operators, and bounded operators in X are denoted by $\tilde{C}(X)$, $C(X)$, $B(X)$, respectively.  

Let $A$ be a linear relation in the Krein space $\mathcal{K}$. When $X=Y=\mathcal{K}$, we use the notation $A^{+}$ rather than $A^{\ast}$. We say that $A$ is \textit{symmetric} (\textit{selfadjoint}) if it satisfies $A\subseteq A^{+}$ ($A=A^{+})$.  

Every point $\alpha \in \mathbb{C}$ for which $\left\{ f,\alpha f \right\}\in A$, with some $f\ne 0$, is called a \textit{finite eigenvalue}, denoted by $\alpha \in \sigma_{p}(A)$, and the vector $f$ is called \textit{eigenvector at the eigenvalue} $\alpha $. If for some $z \in \mathbb{C}$, the operator $(A-z)^{-1}$ is bounded, not necessarily densely defined in $\mathcal{K}$, then $z$ is a \textit{point of regular type of} $A$, symbolically, $z\in \hat{\rho }\left( A \right)$. If for $z\in \mathbb{C}$, the relation $\left( A-z \right)^{-1}$ is a bounded operator and $\overline{\ran \left(A-z \right)} = \mathcal{K}$ ($\ran \left(A-z \right)= \mathcal{K}$, if $A$ is closed), then $z$ is a \textit{regular point of} $A$, symbolically $z\in \rho \left( A \right)$.  
\\

\textbf{1.3. Introduction}

Let $(\mathcal{K}\left[ .,. \right])$ be a separable Krein (Pontryagin or Hilbert) space. Typically, the concept of a boundary triple $\Pi=(\mathcal{H},\Gamma_{0}, \Gamma_{1})$ is introduced by assuming the existence of a closed symmetric linear operator or relation $S\subseteq \mathcal{K}^{2}$, with equal, finite, or infinite defect numbers, c.f. \cite{DM1,BeLu,BHS,DHMS1,DHMS2,BDHS}. Over time, the concept evolved from the ordinary boundary triple $\Pi$ for a symmetric operator $S$ in a Hilbert space $\mathcal{K}$ with a surjective reduction operator $\Gamma$, to generalized versions with a symmetric relations $S$ in a Krein space $\mathcal{K}$ with reduction relation $\Gamma$; see e.g., \cite{DM2, BHS, DHM}. In all above papers, the domain of the reduction operator or relation $\Gamma=\lbrace\Gamma_{0}, \Gamma_{1}\rbrace : \mathcal{K}^{2} \rightarrow \mathcal{H}^{2}$ was defined by means of $S$, precisely by means of the adjoint relation $S^{+}\subseteq \mathcal{K}^{2}$. The following is one of those definitions. 

\begin{definition}\label{definition12} \cite{D1} Let $S$ be a closed Symmetric relation on a seperable Krein space. A triple $\Pi =(\mathcal{H}, \Gamma _{0}, \Gamma_{1})$, where $\mathcal{H}$ is a Hilbert space and $\Gamma_{0}, \Gamma_{1}$ are bounded operators from $S^{+}$ to $\mathcal{H}$, is called an ordinary boundary triple for the relation $S^{+}$ if the abstract Green's identity 
\begin{equation}
\label{eq12}
\left[ f',g \right]-\left[ f,g' \right]=\left( \Gamma_{1}\hat{f},\, 
\Gamma_{0}\hat{g} \right)_{\mathcal{H}}-\left( \Gamma_{0}\hat{f},\, \Gamma 
_{1}\hat{g} \right)_{\mathcal{H}}, \forall \hat{f}, \hat{g}\in S^{+}, 
\end{equation}
holds, and the mapping $\Gamma :\hat{f}\to \left( {\begin{array}{*{20}c}
\Gamma_{0}\hat{f}\\
\Gamma_{1}\hat{f}\\
\end{array} } \right)$ from $S^{+}$ to $\mathcal{H}\times \mathcal{H}$ is surjective. The operator $\Gamma $ is called boundary or reduction operator.
\end{definition}

In other kinds of boundary triples $\Pi =(\mathcal{H}, \Gamma _{0}, \Gamma_{1})$, the existence of a symmetric relation $S$ is also assumed; see \cite{BL,DHM}. However, in applications we usually do not know the symmetric relation or operator $S$ in advance. We usually know only a differential expression and the corresponding differential equation. 

One of the simplest examples is as follows. Let us denote $\mathcal{K}=L^{2}(0,\infty )$. Then the equation is 
\begin{equation}
\label{eq14}
\frac{d^{2}y}{dx^{2}}-\lambda y=g\left( x \right),
\end{equation}
and the differential expression is $l:=-\frac{d^{2}}{dx^{2}}$, $x \in \left( 0,\infty \right)$. The expression $l$ is regular at $0$ and the limit-point case at $\infty $. We know that the maximal domain $\mathcal{D}$ of $l$ consists of all functions $f\in L^{2}(0,\infty )$ for which $f$ and $\frac{df}{dx}$ are absolutely continuous, and $\frac{d^{2}f}{dx^{2}}\in L^{2}(0,\infty )$. We do not know in advance what might be the closed linear relation or operator $S\subseteq \mathcal{K}^{2}$.  Then, if boundary operator $\Gamma =\lbrace \Gamma _{0}, \Gamma_{1}\rbrace$ is defined by
\begin{equation}
\label{eq16}
\mathcal{K}:=L^{2}(0;\infty ); \mathcal{H}=\mathbb{C};  \Gamma_{0}f:=f\left( 0 \right); \Gamma_{1}f:=f'\left( 0 \right),
\end{equation}
one finds out that $S:=\lbrace f \in \mathcal{D} : f\left( 0 \right)=0 \wedge f'\left( 0 \right)=0 \rbrace$ is the closed symmetric linear relation that satisfies the abstract Green's identity (\ref{eq12}). 

Our point is that in applications, such as this example, we do not have $S$ in advance, whereas in theory, we assume the existence of $S$ in advance. We can even determine the Weyl function and the Weyl solution of the equation $\frac{d^{2}y}{dx^{2}}-\lambda y=0$ and the defect number without ever determining the symmetric linear relation $S$. 

This motivates us to introduce a general model using a relation $\Gamma: \mathcal{K}^{2} \rightarrow \mathcal{H}^{2}$ that satisfies the abstract Green's identity, without assuming in advance the existence of the closed symmetric linear relation $S$. 

In \cite[Theorem 2.1.(b)]{B2}, we began with a strict and regular function $Q\in N_{\kappa }\left( \mathcal{H} \right)$ and proved the existence of a closed symmetric relation $S$ with a corresponding ordinary boundary triple $\Pi=(\mathcal{H}, \Gamma_{0}, \Gamma_{1})$ such that $Q$ was the Weyl function associated with $S$ and $\Pi$. The key point is that again we did not have to assume the existence of the linear relation $S$ in order to obtain an ordinary boundary triple.

Motivated by the above two cases, we introduce the Green's boundary relation, abbreviated by GBR, in Definition \ref{definition24}, where we only assume the existence of a separable Krein space $\left( \mathcal{K}, \left[ .,. \right] \right)$, the existence of an auxiliary Hilbert space $\left( \mathcal{H}, \left( .,. \right) \right)$, and the existence of a linear relation $\Gamma : \mathcal{K}^{2} \rightarrow\mathcal{H}^{2}$ that satisfies the abstract Green's identity (\ref{eq212}).

In Section \ref{s4}, in addition to Definition \ref{definition24}, some properties of the Green's boundary relation are proved. For example, in Proposition \ref{proposition28}, a condition on the GBR that is equivalent to the existence of the symmetric closed relation $S\subseteq \mathcal{K}^{2}$ that satisfies $ S^{+}=\overline{\dom}\, \Gamma \subseteq \mathcal{K}^{2}$ is given. This means that a general model is created, the Green's boundary model, in which we can fit all kinds of already defined boundary triples, including ordinary, isometric, unitary and others; see e.g., \cite{DM1, DM2, DHM}. 

In Section \ref{s6}, we utilize the Krein space $\left( X\times Y, \left[ .,. \right]_{X\times Y} \right)$, where $X$ and $Y$ are Krein spaces, to establish a couple of generalizations from Hilbert spaces to Krein spaces. We prove Proposition \ref{proposition34}, which generalizes \cite[Proposition 1.3.2]{BDHS}. We also prove Propositions \ref{proposition310}, a generalization of \cite[Proposition 2.7 (i)]{DHMS1}.

In Section \ref{s8}, in Proposition \ref{proposition42}, some general properties of the isometric and unitary linear relation $V:X\to Y$, where $X$ and $Y$ are Krein spaces, that we need in the research of GBR, are proved. For example, Proposition \ref{proposition42} (ii) is generalization of \cite[Lemma 1.8.1]{BHS}. Similar  statements to Proposition \ref{proposition42}, but slightly more special, can be found in \cite[Corollary 2.4]{DHMS1}. 

In Section \ref{s10}, Proposition \ref{proposition52}, the results of Section \ref{s8} and properties of the space $\left( X\times Y, \left[ .,. \right]_{X\times Y} \right)$, when $X=\mathcal{K}^{2}, Y=\mathcal{H}^{2}, V=\Gamma$, are applied on the Green's boundary relation $\Gamma : \mathcal{K}^{2} \rightarrow \mathcal{H}^{2}$ in order to prove some important properties of isometric and unitary boundary relations $\Gamma $. 

If $\Gamma: \mathcal{K}^{2} \rightarrow \mathcal{H}^{2} $ is a reduction operator of an ordinary boundary triple $\Pi$, then the closed symmetric linear relation $S$ associated with $\Pi$ coincides with $N:= \ker \Gamma $ and with the isotropic manifold $M$ of $\dom\, \Gamma $. However, this is not necessarily the case for more general Green's boundary relations. In Corollary \ref{corollary56}, we provide sufficient conditions for $S=N=M$ if $\Gamma$ is a Green's boundary relation. Additionally, in Proposition \ref{proposition510}, we prove $S=N=M$ if $\Gamma$ is a unitary boundary relation. 

In subsection 5.2, we utilize the well known main transformation $\mathcal{J}$; please refer to \cite{DHMS1,BDHS} for defintion and basic propertis of $\mathcal{J}$. For the convenience of the reader, we have repeated the definition of $\mathcal{J}$, see (\ref{eq214}). We employ the transformation $\mathcal{J}$ to further study isometric and unitary boundary relations, denoted by $\Gamma$, as well as symmetric and self-adjoint relations, $\tilde{A}=\mathcal{J}(\Gamma)$. For instance, we establish a surprising result showing that $\tilde{A}=\mathcal{J}(\Gamma)$ is reduced through the decomposition $ \tilde{\mathcal{K}}=\mathcal{K} [+] \mathcal{H}$ if and only if $\Gamma = \ker\, \Gamma \times \mul\, \Gamma $, that is, if and only if $\Gamma $ is a \textit{trivial} unitary relation. 

Through a couple of examples that follow, we investigate various scenarios distinct from the one presented in the preceding result. In doing so, we aim to analyze the role of specific assumptions underlying the stated result. In Example \ref{example520}, we utilize the differential equation provided by (\ref{eq14}) and the triple specified by (\ref{eq16}) to construct an irreducible self-adjoint operator $\tilde{A}$. Additionally, we examine a scenario where the vector spaces $\mathcal{K}$ and $\mathcal{H}$ satisfy $\mathcal{K} = \mathcal{H}$, thus diverging from the standard condition $\mathcal{K} \cap \mathcal{H} = {0}$. 

In Section \ref{s12}, the previously proven properties of Green's boundary relation are used to generalize two statements about generalized Nevanlinna families. One generalization is Proposition \ref{proposition64} and the other is Theorem \ref{theorem66} which generalizes \cite[Theorem 4.8]{BDHS}. 

In Section \ref{s14}, we study Green's boundary relations $\Gamma$ with non-degenerate $\ran\, \Gamma$. For example, we study the closure $\bar{\Gamma}$ of the Green's boundary relation with non-degenerate  $\ran\, \Gamma$. We also provide conditions that make the Green's boundary relation becomes an ordinary boundary relation.

In Section \ref{s16}, we show another use of the Green's boundary relation. The well-known boundary triples, previously studied in \cite{DHM}, are fitted into Green's boundary relation. That is, the boundary triples with Hilbert space $\mathcal{K}$, such as ordinary, isometric, unitary, AB-generalized, B-generalized, quasi-boundary and S-boundary relations are generalized in terms of the Green's boundary relation to boundary triples with Krein space $\mathcal{K}$.
\section{Green's boundary model}\label{s4}

\textbf{2.1.} Let $(\mathcal{K},\left[ .,. \right])$ and $(\mathcal{H},(.,.))$ be Krein and Hilbert spaces, respectively, and let $(\mathcal{K}^{2}, \left[ .,. \right])\,$, and $(\mathcal{H}^{2}, (.,.))$ be the corresponding \textit{product spaces}, with inner products $\left[ .,. \right]$ and $(.,.)$ defined as follows:
\begin{equation}
\label{eq22}
\left[ \hat{f},\hat{g} \right]=\left[ f,g \right]+\left[ f',g' 
\right],\hat{f}=\left( {\begin{array}{*{20}c}
f\\
f'\\
\end{array} } \right)\in \mathcal{K}^{2}, \, \hat{g}=\left( {\begin{array}{*{20}c}
g\\
g'\\
\end{array} } \right)\in \mathcal{K}^{2}, 
\end{equation}
\begin{equation}
\label{eq24}
\left( \hat{h},\hat{k} \right)=\left( h,k \right)+\left( h',k' 
\right), \hat{h}=\left( {\begin{array}{*{20}c}
h\\
h'\\
\end{array} } \right)\in \mathcal{H}^{2}, \, \hat{k}=\left( {\begin{array}{*{20}c}
k\\
k'\\
\end{array} } \right)\in \mathcal{H}^{2}.
\end{equation}
Let us introduce the involution
\begin{equation}
\label{eq26}
J_{\mathcal{K}}:=\left( {\begin{array}{*{20}c}
0 & -iI_{\mathcal{K}}\\
iI_{\mathcal{K}} & 0\\
\end{array} } \right), \, J_{\mathcal{H}}:=\left( {\begin{array}{*{20}c}
0 & -iI_{\mathcal{H}}\\
iI_{\mathcal{H}} & 0\\
\end{array} } \right),
\end{equation}
where $I_{\mathcal{K}}$ and $I_{\mathcal{H}}$ are identities in the corresponding spaces, and let $\left( \mathcal{K}^{2},\left[ .,. \right]_{\mathcal{K}^{2}} \right)$ and $\left( \mathcal{H}^{2}, \left[ .,. \right]_{\mathcal{H}^{2}} \right)$ be the inner product spaces obtained by means of the inner products 
\begin{equation}
\label{eq28}
\left[ \hat{f},\hat{g} \right]_{\mathcal{K}^{2}}:=\left[ J_{\mathcal{K}}\hat{f},\hat{g} 
\right];\, \hat{f},\hat{g}\in \mathcal{K}^{2},
\end{equation}
\begin{equation}
\label{eq210}
\left[ \hat{h},\hat{k} \right]_{\mathcal{H}^{2}}:=\left( J_{\mathcal{H}}\hat{h},\hat{k} 
\right); \,  \hat{h},\hat{k}\in \mathcal{H}^{2},
\end{equation} respectively. Frequently, the notations $\mathcal{K}^{2}=\left( \mathcal{K}^{2}, \left[ .,. \right]_{\mathcal{K}^{2}}\right)=\left( \mathcal{K}^{2},J_{\mathcal{K}}\right) $ and $\mathcal{H}^{2}=\left( \mathcal{H}^{2}, \left[ .,. \right]_{\mathcal{H}^{2}}\right)=\left( \mathcal{H}^{2},J_{\mathcal{H}}\right) $ are used. 

In the following lemma we collect several well known facts.

\begin{lemma}\label{lemma22} Let $(\mathcal{K},\left[ .,. \right])$ be a Krein space and let $\left( \mathcal{K},( .,.) \right) $ be the Hilbert space associated with $(\mathcal{K},\left[ .,. \right])$ by means of the fundamental symmetry $G$, i.e. $\left[x,y \right]=\left( Gx,y \right),  \forall x,y \in \mathcal{K}$. The spaces $\left( \mathcal{K}^{2}, \left[ .,. \right]_{\mathcal{K}^{2}}\right) $ and $\left( \mathcal{H}^{2},  \left[ .,. \right]_{\mathcal{H}^{2}}\right) $ with scalar products (\ref{eq28}) and (\ref{eq210}), respectively, are Krein spaces. They are associated with the product Hilbert spaces $(\mathcal{K}^{2}, \left( .,. \right))$ and $(\mathcal{H}^{2},\left( .,. \right))$ with the fundamental symmetries
\[
\mathcal{G}=\left( {\begin{array}{*{20}c}
0 & -iG\\
iG & 0\\
\end{array} } \right) \,\wedge \, \mathcal{J}=\left( {\begin{array}{*{20}c}
0 & -iI_{\mathcal{H}}\\
iI_{\mathcal{H}} & 0\\
\end{array} } \right)
\]
respectively. Symbolically,
\begin{equation}
\label{eq211.5}
\left[ \hat{f},\hat{g} \right]_{\mathcal{K}^{2}}=\left(\mathcal{G} \hat{f},\hat{g} \right)_{\mathcal{K}^{2}}, \forall \hat{f}, \hat{g} \in \mathcal{K}^{2}.
\end{equation}
and (\ref{eq210}) hold.
\end{lemma}
If $\left( X,  {[.,.]}_{X} \right)$ and $\left(Y, {[.,.]}_{Y} \right)$ are Krein spaces, a linear relation $V\subseteq X\times Y$ is called \textit{isometric} if $V^{-1}\subseteq V^{\ast }$; see \cite[p.5360]{DHMS1}.  Obviously, it is equivalent with 
\[
\left[ f' ,g' \right]_{Y}=\left[ f,g \right]_{X},\forall 
\left( {\begin{array}{*{20}c}
f\\
f'\\
\end{array} } \right),\, \left( {\begin{array}{*{20}c}
g\\
g'\\
\end{array} } \right)\in V.
\]
A linear relation $V\subseteq X\times Y$ is called \textit{unitary} if $V^{-1}=V^{\ast }$. 

\begin{definition}\label{definition24} Let $\Gamma :\hat{f}\to \hat{h}, \hat{f} \in \mathcal{K}^{2},\hat{h} \in \mathcal{H}^{2}$, be a linear relation (operator) from the Krein space $(\mathcal{K}^{2},\left[ .,. \right]_{\mathcal{K}^{2}})$ to Krein space $\left( \mathcal{H}^{2},\left[ .,. \right]_{\mathcal{H}^{2}} \right)$ which satisfies 
\begin{equation}
\label{eq212}
\left[ f',g \right]-\left[ f,g' \right]=\left( h',k \right)-\left( 
h,k' \right),\forall \left( {\begin{array}{*{20}c}
\hat{f}\\
\hat{h}\\
\end{array} } \right), \left( {\begin{array}{*{20}c}
\hat{g}\\
\hat{k}\\
\end{array} } \right)\in \Gamma \subseteq \mathcal{K}^{2}\times H^{2},
\end{equation}
or equivalently
\begin{equation}
\label{eq214}
\left[ \hat{f},\hat{g} \right]_{\mathcal{K}^{2}}=\left[ \hat{h},\hat{k} 
\right]_{\mathcal{H}^{2}},\forall \left( {\begin{array}{*{20}c}
\hat{f}\\
\hat{h}\\
\end{array} } \right),\left( {\begin{array}{*{20}c}
\hat{g}\\
\hat{k}\\
\end{array} } \right)\in \Gamma \subseteq \mathcal{K}^{2}\times \mathcal{H}^{2}.
\end{equation}
Then the relation $\Gamma $ is called \textbf{Green's boundary relation (operator)}.
\end{definition}
Note that the Green's boundary relation, abbreviated as GBR, $\Gamma :\mathcal{K}^{2}\to \mathcal{H}^{2}$ introduced here is a \textbf{generalization} of the boundary relations defined in \cite[Definition 3.1]{DHMS1} and \cite[Definition 3.1]{BDHS} in the following aspects:
\begin{itemize}
\item We do not have to assume the existence of a closed symmetric linear relation $S\subseteq \mathcal{K}^{2}$ in order to define a GBR.
\item We do not set any specific assumptions either on $\dom\, \Gamma \subseteq \mathcal{K}^{2}$ or on $\ran\, \Gamma \subseteq \mathcal{H}^{2}$. 
\item We do not assume that the condition \cite[Definition 3.1.G2]{DHMS1} is satisfied, i.e. we do \textbf{not} assume that $\Gamma $ is a unitary relation.  
\end{itemize}

The Green's boundary relation $\Gamma$ is associated with two relations defined by
\begin{equation}
\label{eq220}
\Gamma_{0}=\left\{ \left( {\begin{array}{*{20}c}
\hat{f}\\
h\\
\end{array} } \right):\left( {\begin{array}{*{20}c}
\hat{f}\\
\hat{h}\\
\end{array} } \right)\in  \Gamma  \right\}\, \wedge \, \Gamma_{1}=\left\{ \left( {\begin{array}{*{20}c}
\hat{f}\\
h'\\
\end{array} } \right):\, \left( {\begin{array}{*{20}c}
\hat{f}\\
\hat{h}\\
\end{array} },  \right)\in \, \Gamma \right\}.
\end{equation} 

The setup described in Section \ref{s4}.1 will be referred to as the \textit{general Green's boundary model}, or simply the \textit{Green's boundary model}. It is called "general" because we did not define it using a symmetric closed linear relation $S \subseteq \mathcal{K} \times \mathcal{K}$ as is typically done. 
\\

\textbf{2.2.} The following proposition is probably known for all boundary triples defined by means of a closed symmetric linear relation $S$. Here, we establish it without the assumption that $S$ exists.

\begin{proposition}\label{proposition26} Let $\Gamma :\mathcal{K}^{2} \to \mathcal{H}^{2}$ be a Green's boundary relation, then:
\begin{enumerate}[(i)]%, (i), (ii),...
\item $\ker \Gamma_{i} \subseteq {\left( \ker \Gamma_{i} \right)}^{+}, i=0, 1$,
\item $\ker \Gamma \subseteq \left( \ker \Gamma \right)^{+}$,
\item $\mul\, \Gamma_{i} \subseteq \left(\mul\,\Gamma_{i} \right)^{+}, i=0, 1$,
\item $\mul\, \Gamma \subseteq {\left(\mul\,\Gamma \right)}^{+}$.
\end{enumerate}
\end{proposition}
\textbf{Proof.}

(i) If $\hat{f}, \hat{g}\in \ker \Gamma_{0}$, then 
$\left\lbrace \hat{f} , \left( {\begin{array}{*{20}c}
0\\
h'\\
\end{array} } \right)\right\rbrace  \in \Gamma $ and $\left\lbrace \hat{g} , \left( 
{\begin{array}{*{20}c}
0\\
k'\\
\end{array} } \right)\right\rbrace  \in \Gamma$. Assume that $\hat{f}\in \ker \Gamma_{0}$ is 
arbitrarily selected. Then for every $\hat{g}\in \ker \Gamma_{0}$, 
according to (\ref{eq212}), we have 
\[
\left( \left[ f',g \right]-\left[ f,g' \right]=0 \right)\Rightarrow \hat{f}\in \left( \ker \Gamma_{0} \right)^{+}.
\]
Hence, the relation $\ker \Gamma_{0}$ is symmetric. Similar reasoning can be repeated for $\ker \Gamma_{1}$.

Claims (ii), (iii) and (iv) follow also from (\ref{eq212}) directly. \hfill$\square$ 
\\

In the following proposition, we give a connection between Green's boundary relation and boundary relations defined by means of the closed symmetric linear relation $S$ in earlier papers. 

\begin{proposition}\label{proposition28} Let $\Gamma :\mathcal{K}^{2} \to \mathcal{H}^{2}$ be a Green's boundary relation. Then:

\begin{enumerate}[(i)]%, (i), (ii),...
\item The condition 
\begin{equation}
\label{eq222}
\left( \dom\, \Gamma \right)^{\left[ \bot \right]_{\mathcal{K}^{2}}} \subseteq \overline{\dom}\, \Gamma 
\end{equation}
or, equivalently,
\[
\left( \left[ f',g \right]-\left[ f,g' \right]=0,\forall \hat{g}\in \dom\, \Gamma \right)\Rightarrow \hat{f}\in  \overline{\dom}\, \Gamma
\]
holds if and only if there exists a closed symmetric linear relation 
$S\subseteq \mathcal{K}^{2}$ such that $S^{+}=\overline{\dom}\, \Gamma $.
\item If 
\begin{equation}
\label{eq223}
\left( \dom\, \Gamma \right)^{\left[ \bot \right]_{\mathcal{K}^{2}}}\subseteq \dom\, \Gamma ,
\end{equation}
then $M=S$, where $M$ denotes the isotropic manifold of $dom\, \Gamma$.
\end{enumerate}
\end{proposition}
\textbf{Proof.} 

(i) Assume that $\Gamma :\mathcal{K}^{2} \to \mathcal{H}^{2}$ satisfies (\ref{eq222}). We define 
\begin{equation}
\label{eq224}
S:=\left(\dom\, \Gamma  \right)^{+}=\left( \dom\, \Gamma \right)^{\left[ \bot \right]_{\mathcal{K}^{2}}}.
\end{equation}
Hence, $S$ is a closed relation. It follows 
\[{S}^{+}=\left(\left( \dom\, \Gamma \right)^{\left[ \bot \right]_{\mathcal{K}^{2}}}\right)^{\left[ \bot \right]_{\mathcal{K}^{2}}} = \overline{\dom}\, \Gamma.
\] It remains to see $S\subseteq {S\, }^{+}$. According to definition of $S$ and (\ref{eq222})
\[
\hat{g} \in S \Rightarrow \hat{g} \in \left( \dom\, \Gamma \right)^{\left[ \bot \right]_{\mathcal{K}^{2}}} \Rightarrow \hat{g} \in  \overline{\dom}\, \Gamma=S^{+}.
\]

Conversely, assume that $S$ is a closed,  $S\subseteq S^{+}$ and $S^{+}=\overline{\dom}\, \Gamma $. Then obviously $S=\left(\overline{\dom}\, \Gamma \right)^{+}$. To prove $(\ref{eq222})$, assume 
\[
\left[ f',g \right]-\left[ f,g' \right]=\left[ \hat{f},\hat{g} \right]_{\mathcal{K}^{2}}=0,\forall \hat{g}\in \dom\, \Gamma .
\]
According to Lemma \ref{lemma22}, the scalar product $\left[ .,. \right]_{\mathcal{K}^{2}}$ is continuous with 
respect to the strong topology of the product Hilbert space $\left( \mathcal{K}^{2}, \left( .,. \right) \right)$. Thus $\left[ \hat{f},\hat{g} \right]_{\mathcal{K}^{2}}=0,\, \forall \hat{g}\in \overline{\dom}\,\Gamma $. This further means 
\[
\hat{f}\in \left( \overline{\dom}\, \Gamma  \right)^{+}=S \subseteq S^{+}=\overline{\dom}\, \Gamma.
\]
This proves implication (\ref{eq222}).

(i) By definition 
\[
M=\dom\, \Gamma \cap \left(\dom\, \Gamma \right)^{\left[ \bot \right]_{\mathcal{K}^{2}}}.
\]
From this, (\ref{eq223}) and (\ref{eq224}) we get 
\[
M=\left(\dom\, \Gamma \right)^{\left[ \bot \right]_{\mathcal{K}^{2}}}=S .
\]\hfill$\square$ 

Note that the condition (\ref{eq222}) determines is $S$ an operator or a multi-valued relation; we do not have to guess in advance which one to assume. 

\begin{definition}\label{definition210} Green's boundary relation $\Gamma :\mathcal{K}^{2} \to \mathcal{H}^{2}$ that satisfies the maximality condition (\ref{eq222}) is called an isometric \textbf{boundary} relation. Then the relation $S$ defined by (\ref{eq224}) is called a \textit{closed symmetric relation associated with} $\Gamma$.
\end{definition}

Note that this definition is equivalent with usual definition of an isometric boundary relation. We rewrote it here to place it in the context of Green's boundary model. Also note that every Green's boundary relation $\Gamma$ is an isometric relation between Krein spaces $\mathcal{K}^{2} $ and $\mathcal{H}^{2}$, but it is not an isometric \textbf{boundary} relation if it does not satisfy (\ref{eq222}).
 
\begin{definition}\label{definition214} The Green's boundary relation $ \Gamma :\mathcal{K}^{2} \to \mathcal{H}^{2}$  that satisfies:
\[
\Gamma^{-1}=\Gamma^{\ast} 
\]
is called \textit{unitary boundary relation}.
\end{definition}
Later, in Proposition \ref{proposition510} we will prove that $\Gamma ^{-1}=\Gamma ^{\ast }$ implies (\ref{eq223}), which means that every unitary boundary relation must have a corresponding symmetric relation $S$. Therefore, this definition is equivalent to the usual definition of an unitary boundary relation even though here we do not assume existence of the relation $S$ in advance.

Proposition \ref{proposition28} (i) shows us that the condition (\ref{eq222}) is essential because it makes the relation $S$ defined by (\ref{eq224}) a symmetric relation. Proposition \ref{proposition28} (ii) gives us a sufficient condition for $S=M$. Proposition \ref{proposition52} (i) will give us a sufficient condition for $M=N:=\ker\, \Gamma$. Corollary \ref{corollary56} will give us a couple of sufficient condition for $S=N=M$.

\begin{corollary}\label{corollary216} Let $\Gamma : \mathcal{K}^{2} \to \mathcal{H}^{2}$ be a Green's boundary relation. \begin{enumerate}[(i)]%, (i), (ii),...
\item If there exists a symmetric extensions $A$ of $S:={\left( \dom\,\Gamma \right)}^{+}$, then $\Gamma$ satisfies (\ref{eq222}), i.e. $\Gamma$ is an isometric boundary relation. 
\item If additionally $A^{+} \subseteq \dom \, \Gamma$, then $\tilde{\Gamma }:=\Gamma_{\mid A^{+}} $ satisfies (\ref{eq223}). 
\end{enumerate}
\end{corollary}
\textbf{Proof.} 

(i) We assume $S={\left( \dom\, \Gamma \right)}^{+} \subseteq A\subseteq A^{+}$ and $A^{+}\subseteq S^{+}=\overline{\dom}\, \Gamma $. Hence, $S$ is also symmetric. According to Proposition \ref{proposition28} (i), $\Gamma$ satisfies (\ref{eq222}).

(ii) From $\dom\, \tilde{\Gamma}=A^{+}$ it follows
\[
\hat{f} {\left[ \bot \right]_{\mathcal{K}^{2}}} \dom\, \tilde{\Gamma}=A^{+} \Rightarrow \hat{f}\in {A^{+}}^{+}=\bar{A}\subseteq A^{+} =\dom\, \tilde{\Gamma} .
\]
\hfill  $\square$ 
\section{Krein space $\left( X\times Y, \left[ .,. \right]_{X\times Y} \right)$}\label{s6}
\textbf{3.1.} Let $\left( X, \left( .,. \right)_{X} \right)$ and $\left( Y, \left( .,. \right)_{Y} \right)$, $X \cap Y =\lbrace 0 \rbrace$, be the Hilbert spaces associated with the Krein spaces $\left( X, {[.,.]}_{X} \right)$ and $\left( Y, {[.,.]}_{Y} \right)$, respectively. This means 
\[
\left[ u_{x},v_{x} \right]_{X}=\left( J_{X}u_{x},v_{x} \right)_{X},\, \left[ u_{y},v_{y} \right]_{Y}=\left( J_{Y}u_{y},v_{y} \right)_{Y}; \, \left( {\begin{array}{*{20}c}
u_{x}\\
u_{y}\\
\end{array} } \right), \left( {\begin{array}{*{20}c}
v_{x}\\
v_{y}\\
\end{array} } \right)\in X\times Y,
\]
where $J_{X}$ and $J_{Y}$ are fundamental symmetries of the respective Krein spaces $X$ and $Y$. The \textit{product Hilbert space} $\left( X\times Y,\left( .,. \right)_{X\times Y} \right)$ is defined as usual
\begin{equation}
\label{eq31}
\left( \left( {\begin{array}{*{20}c}
u_{x}\\
u_{y}\\
\end{array} } \right),\left( {\begin{array}{*{20}c}
v_{x}\\
v_{y}\\
\end{array} } \right) \right)_{X\times Y}:=\left( u_{x},v_{x} 
\right)_{X}+\left( u_{y},v_{y} \right)_{Y};\, \, \left( 
{\begin{array}{*{20}c}
u_{x}\\
u_{y}\\
\end{array} } \right),\left( {\begin{array}{*{20}c}
v_{x}\\
v_{y}\\
\end{array} } \right)\in X\times Y.
\end{equation}
We will need the Krien space $\left( X\times Y, \left[ .,. \right]_{X\times Y} \right)$. It was introduced by \cite[(1)]{Sh}, i.e.  
\begin{equation}
\label{eq32}
\left[\left( {\begin{array}{*{20}c}
u_{x}\\
u_{y}\\
\end{array} } \right),\left( {\begin{array}{*{20}c}
v_{x}\\
v_{y}\\
\end{array} } \right) \right]_{X\times Y}:=\left[ u_{x},v_{x} 
\right]_{X}-\left[ u_{y},v_{y} \right]_{Y}; \left({\begin{array}{*{20}c}
u_{x}\\
u_{y}\\
\end{array} } \right),\left( {\begin{array}{*{20}c}
v_{x}\\
v_{y}\\
\end{array} } \right)\in X\times Y.
\end{equation} 

\begin{lemma}\label{lemma32} The space $\left( X\times Y,\left[ .,. \right]_{X\times Y} \right)$, defined by (\ref{eq32}), is a Krein space with the fundamental symmetry 
\[
J_{X \times Y}:=\left( {\begin{array}{*{20}c}
J_{X} & 0\\
0 & -J_{Y}\\
\end{array} } \right)
\]
with respect to the \textbf{product} Hilbert space $\left( X\times Y,{(.,.)}_{X\times Y} \right)$ defined by (\ref{eq31}). 
\end{lemma}
\textbf{Proof.} It is easy to verify $\left[ \left( {\begin{array}{*{20}c}
u_{x}\\
0\\
\end{array} } \right),\left( {\begin{array}{*{20}c}
0\\
v_{y}\\
\end{array} } \right) \right]_{X\times Y}=0$. Therefore, we can consider 
\[
X\times Y=X\left[ + \right]Y,
\]
where $\left[ + \right]$ is the direct and orthogonal sum with respect to scalar product (\ref{eq32}). The space $X\times Y$ is a Krein space as an orthogonal sum of Krein spaces. \hfill$\square$

\begin{remark}\label{remark33} Let us note that $\left( X\times Y,{[.,.]}_{X\times Y} \right)$ is still a Krein space, even if Krein spaces $X$ and $Y$ are equal as vector spaces, but with different scalar products. However, some properties of the Krein space $\left( X\times Y,{[.,.]}_{X\times Y} \right)$ with $X=Y$ may be different from the properties of the Krein space where vector spaces $X$ and $Y$ satisfy $X \cap Y=\lbrace 0\rbrace $; see Lemma \ref{lemma518} and Example \ref{example519}.
\end{remark}

\textbf{3.2.} The following proposition for a linear relation $H\subseteq X\times Y$ can be found for Hilbert spaces in \cite[Proposition 1.3.2]{BHS}. Also, all statements of the following proposition are known when Krein spaces $X$ and $Y$ coincide. Statements (iv) and (v) of the following proposition have been proven for any Krein spaces $X$ and $Y$ in \cite[(2.7)]{DHMS1} by means of the corresponding statements for Hilbert spaces and fundamental symmetries. We will prove all the statements of the following proposition directly in terms of the Krein spaces $X$ and $Y$, by means of the Krein space $\left( X\times Y, \left[ .,. \right]_{X\times Y} \right)$.

As usually $\bar{H}$ denotes the closure of $H\subseteq X\times Y$ with respect to product Hilbert space $X \times Y$. 

\begin{proposition}\label{proposition34} Let $\left( X, {[.,.]}_{X} \right)$ 
and $\left( Y, {[.,.]}_{Y} \right)$ be Krein spaces, and let $H:X\to Y$ 
be a linear relation. Let $H^{\ast }:Y \to X$ be the adjoint linear 
\textbf{relation} of $H$. Then 

\begin{enumerate}[(i)]%, (i), (ii),...
\item $H^{\ast }$ is a closed linear relation,
\item $\left( \bar{H} \right)^{\ast }=H^{\ast }$, 
\item $\bar{H}={H^{\ast }}^{\ast }$,
\item $\left( \dom\, H\, \right)^{\left[ \bot \right]_{X}}=\mul\, H^{\ast }$, 
\item $\ker\, H^{\ast }=\left( \ran\, H\, \right)^{\left[ \bot \right]_{Y}}$, 
\item $\overline{\dom}\, \bar{H }=\overline{\dom}\, H \, \wedge \, \overline{\ran}\, \bar{H }=\overline{\ran}\, H$.
\end{enumerate}
\end{proposition}
\textbf{Proof}.

(i) From definition of $H^{\ast }$, it follows: 
\begin{equation}
\label{eq36}
H^{\ast }={ \left( H^{\left[ \bot \right]_{X\times Y}} \right)}^{-1}.
\end{equation}
Since the orthogonal companion $H^{\left[ \bot \right]_{X\times Y}}$ is closed, $H^{\ast }$ is also closed. 

(ii) Let us now select any $\left( {\begin{array}{*{20}c}
g_{y}\\
g_{x}\\
\end{array} } \right)\in H^{\ast }$. Because of the continuity of scalar 
product (\ref{eq32}), the relation $
\left[ f_{x},g_{x} \right]_{X}-\left[ f_{y},g_{y} \right]_{Y}=0,\, 
\forall \left( {\begin{array}{*{20}c}
f_{x}\\
f_{y}\\
\end{array} } \right)\in H
$ extends to all $\left( {\begin{array}{*{20}c}
f_{x}\\
f_{y}\\
\end{array} } \right)\in \bar{H}$. Hence, $H^{\ast }\subseteq \left( \bar{H} 
\right)^{\ast }$. This together with the obvious $\left( \bar{H} \right)^{\ast 
}\subseteq H^{\ast }$ gives $\left( \bar{H} 
\right)^{\ast }=H^{\ast }$.

(iii) For a linear relation $T\subseteq Y\times X$ one can easily verify
\begin{equation}
\label{eq38}
\left( T^{-1} \right)^{\left[ \bot \right]_{X\times Y}}=\left( T^{\left[ 
\bot \right]_{Y\times X}} \right)^{-1}.
\end{equation}
If we substitute $H^{\ast }\subseteq Y\times X$ for $H \subseteq X \times Y$ into (\ref{eq36}), we obtain
\[
{H^{\ast }}^{\ast }={\, \left( \left( H^{\ast } \right)^{\left[ \bot 
\right]_{Y\times X}} \right)}^{-1}.
\]
If we substitute $H^{\ast }$ for $T$ in (\ref{eq38}) and observe that the 
right sides are equal we get
\[
{H^{\ast }}^{\ast }=\left( \left( H^{\ast } \right)^{-1} \right)^{\left[ \bot 
\right]_{X\times Y}}.
\]
According to (\ref{eq36}) it holds $\left( H^{\ast } \right)^{-1}=H^{\left[ 
\bot \right]_{X\times Y}}$. Finally we get 
\[
{H^{\ast }}^{\ast }=\left( H^{\left[ \bot \right]_{X\times Y}} \right)^{\left[ 
\bot \right]_{X\times Y}}.
\]
Because, $\left( X\times Y,\left[ .,. \right]_{X\times Y} \right)$ is a 
Krein space, all conditions of \cite[p.68, Theorem 6.1]{Bog} are satisfied. Therefore $\left( H^{\left[ \bot \right]_{X\times Y}} \right)^{\left[ \bot \right]_{X\times 
Y}}=\bar{H}$. This proves ${H^{\ast }}^{\ast }=\bar{H}$.

(iv)  
\[
f'\in \left( \dom\, H\, \right)^{\left[ \bot \right]_{X}}\Longleftrightarrow
\left[ h,f' \right]_{X}-\left[ h',0 \right]_{Y}=0,\, \, \forall \left( 
{\begin{array}{*{20}c}
h\\
h'\\
\end{array} } \right)\in H\, \Longleftrightarrow \left( {\begin{array}{*{20}c}
0\\
f'\\
\end{array} } \right)\in H^{\ast}.
\]

(v) $\left( {\begin{array}{*{20}c}
h\\
0\\
\end{array} } \right)\in H^{\ast } \Longleftrightarrow \left[ h,f' 
\right]=0, \forall \left( {\begin{array}{*{20}c}
f\\
f'\\
\end{array} } \right)\in H \Longleftrightarrow h\in \left( \ran\, H\, 
\right)^{\left[ \bot \right]_{Y}}$. 

(vi)  Assume $f\in \dom\, \bar{H}$. Then
$\left( {\begin{array}{*{20}c}
f\\
f'\\
\end{array} } \right)\in \bar{H }$, for some $f'\in \ran\, \bar{H}$. There exists the sequences $\left( \left\{ f_{n} \right\} \subseteq \dom\, H \wedge  \left\{ f_{n}' \right\} \subseteq \ran\, H \right)  $ such that $ \left( f_{n}\to f  \wedge f_{n}' \to f' \right) $, when $n\rightarrow \infty$. This further means $f \in  \overline {\dom}\, H $. Hence, $\dom\, \bar{H } \subseteq \overline{\dom}\, H $. This implies $\overline{\dom}\, \bar{H } \subseteq \overline{\dom}\, H $. Because the converse inclusion is obvious, the identity $\overline{\dom}\, \bar{H } = \overline{\dom}\, H $ holds.

By the same token the identity $\overline{\ran}\, \bar{H } = \overline{\ran}\, H $ holds.  \hfill  $\square $
\\

\textbf{3.3.} The following two propositions are characterizations of isometric and symmetric linear relations in terms of the Krein spaces defined by scalar products (\ref{eq32}) and (\ref{eq28}), respectively.
 
\begin{proposition}\label{proposition38} \cite{Sh} Let $\left( X,\, {[.,.]}_{X} \right)$ and $\left( Y,\, {[.,.]}_{Y} \right)$ be Krein spaces. Then $V\subseteq X\times Y$ is a neutral manifold (hyper-maximal neutral subspace) in the Krein space $\left( X\times Y, \left[ .,. \right]_{X\times Y} \right)$ if and only if $V$ is an isometric (unitary) linear relation. Symbolically,
\[
V \subseteq V^{\left[ \bot\right]_{X \times Y}} \Longleftrightarrow V^{-1} \subseteq V^{*} \left(V=V^{\left[ \bot\right]_{X \times Y}} \Longleftrightarrow V^{-1} = V^{*} \right).
\]
\end{proposition}
%\textbf{Proof}. Assume that $V$ is an isometric relation, i.e. $V^{-1}\subseteq V^{*}$. According to (\ref{eq36}), this is equivalent to $V^{-1} \subseteq { \left( V^{\left[ \bot \right]_{X\times Y}} \right)}^{-1}$. This is further equivalent to $V \subseteq V^{\left[ \bot\right]_{X \times Y}}$. By definition this means that $V\subseteq X\times Y$ is an neutral subspace. 

%The statement for the unitary relation $V$ follows from this if we replace "$\subseteq$" by "$=$" and "neutral manifold" by "hyper-maximal neutral subspace". \hfill $\square$

The following proposition is slight generalization of \cite[Proposition 2.9 (i)]{DHMS1} because $\mathcal{K}$ is a Krein space here rather than a Hilbert space.
  
\begin{proposition}\label{proposition310} Let $(\mathcal{K}, \left[ .,. \right])$ be a Krein space and let $(\mathcal{K}, \left( .,. \right))$ be the corresponding Hilbert space with the fundamental symmetry $G$. Let $(\mathcal{K}^{2}, \left[ .,. \right]_{\mathcal{K}^{2}})$ and $(\mathcal{K}^{2}, \left( .,. \right)_{\mathcal{K}^{2}})$ be spaces as in Lemma \ref{lemma22}. In addition, let $\left( \mathcal{K}\times \mathcal{K}, \left[ .,. \right]_{\mathcal{K}\times \mathcal{K}} \right)$ be the Krein space defined by the scalar product (\ref{eq32}) with $X=Y=\mathcal{K}$.
Then it holds
\begin{enumerate}[(i)]%, (i), (ii),...
\item 
\begin{equation}
\label{eq318}
\left[ \hat{f},\hat{g} \right]_{\mathcal{K}^{2}}= \left[ \left( {\begin{array}{*{20}c}
0 & -i\\
-i & 0\\
\end{array} } \right)\hat{f},\hat{g} \right]_{\mathcal{K}\times \mathcal{K}}, \, \hat{f}, \hat{g} \in\mathcal{K}\times \mathcal{K}.
\end{equation} 
\item A liner relation $A $ in $ (\mathcal{K}, \left[ .,. \right])$ is symmetric if and only if it is a neutral manifold in $(\mathcal{K}^{2}, \left[ .,. \right]_{\mathcal{K}^{2}})$, if and only if $A^{-1} \subseteq A^{\left[ \bot \right]_{\mathcal{K}\times \mathcal{K}}}$.
\item A liner relation $A $ in $ (\mathcal{K}, \left[ .,. \right])$ is self-adjoint if and only if it is a hyper-maximal neutral subspace in $(\mathcal{K}^{2}, \left[ .,. \right]_{\mathcal{K}^{2}})$. Symbolically 
\[
A = A^{+} \Longleftrightarrow A=A^{\left[ \bot \right]_{\mathcal{K}^{2}}} \left( \Longleftrightarrow A^{-1} = A^{\left[ \bot \right]_{\mathcal{K}\times \mathcal{K}}} \right).
\]
\end{enumerate}
\end{proposition}
\textbf{Proof.} 

(i) According to Lemma \ref{lemma32}, the fundamental symmetry of $\mathcal{K} \times \mathcal{K}$ is
\[
J_{\mathcal{K} \times \mathcal{K}}:=\left( {\begin{array}{*{20}c}
G & 0\\
0 & -G\\
\end{array} } \right).
\]
According to Lemma \ref{lemma22}, the fundamental symmetry of $\mathcal{K}^{2}$ is  
\[
\mathcal{G}=\left( {\begin{array}{*{20}c}
0 & -iG\\
iG & 0\\
\end{array} } \right).
\]
The identity (\ref{eq318}) now follows from
\[
\mathcal{G}= J_{\mathcal{K} \times \mathcal{K}}\left( {\begin{array}{*{20}c}
0 & -i\\
-i & 0\\
\end{array} } \right).
\] 

(ii) According to (\ref{eq318}) 
\[
\left({\begin{array}{*{20}c}
f\\
f' \\
\end{array} } \right)  \left[ \bot \right]_{\mathcal{K}^{2}} \left({\begin{array}{*{20}c}
g\\
g' \\
\end{array} } \right) \Longleftrightarrow \left({\begin{array}{*{20}c}
f\\
f' \\
\end{array} } \right)  \left[ \bot \right]_{\mathcal{K}\times \mathcal{K}} \left({\begin{array}{*{20}c}
g'\\
g \\
\end{array} } \right).
\]
Therefore, $A \subseteq A^{+} = A^{\left[ \bot \right]_{\mathcal{K}^{2}}} \Longleftrightarrow A^{-1} \subseteq A^{\left[ \bot \right]_{\mathcal{K}^\times \mathcal{K}}}$.

The statement (iii) follows from (ii) when we substitute "$=$" for "$\subseteq$". \hfill $\square $ 
\\

Note that according to previous two propositions, isometric relations in $\mathcal{K}$ are characterized as neutral in $\left( \mathcal{K} \times \mathcal{K}, \left[ .,.\right]_{\mathcal{K} \times \mathcal{K}} \right) $ and symmetric relations are characterized as neutral in $\left( \mathcal{K}^{2},\left[ .,.\right]_{\mathcal{K}^{2}} \right)$. Therefore, the space $\left( \mathcal{K} \times \mathcal{K}, \left[ .,.\right]_{\mathcal{K} \times \mathcal{K}} \right) $ plays the same role for isometric relations in $\mathcal{K}$ as the space $\left( \mathcal{K}^{2},\left[ .,.\right]_{\mathcal{K}^{2}} \right)$ plays for symmetric relations in $\mathcal{K}$.

\section{Isometric relations between Krein spaces}\label{s8}

The definition of a unitary operator ($V^{-1}=V^{\ast }\, \wedge \mul\, V=\left\{ 0 \right\}$ ) in \cite[p.76]{BHS} is different from the definition in \cite[p.128]{Bog}, where the unitary operator is defined as an isometric operator $V:X\to Y$ with $\dom\, V=X\, \wedge \ran\, V=Y$. The latter operator is a completely invertible isometric, and by many authors it is called \textit{standard unitary operator}. We will say that an isometric (unitary) relation $V \subseteq  X \times Y$ is \textit{trivial} if $V := \ker V \times \mul\, V $.

We will need the statements from the following proposition about the isometric relation $V:X\to Y$, where $X$ and $Y$ are Krein spaces. Recall that we denoted by $M$ the \textit{isotropic manifolds} of $\dom\, V$. Let  $\tilde{M}:= \ran\, V \cap \left( \ran\, V \right)^{[\bot]}$  denote the \textit{isotropic manifolds} of $\ran\, V$. 

\begin{proposition}\label{proposition42} Let $\left( X, {[.,.]}_{X} \right)$ and $\left( Y, {[.,.]}_{Y} \right)$ be Krein spaces, and let $V:X\to Y$ be an isometric relation. Then:
\begin{enumerate}[(i)]%, (i), (ii),...
\item $\tilde{M} = \mul\, V \Leftrightarrow M=\ker V$.
%\item If $\ran\, V$ ($\dom\, V$) does not degenerate, then $M=\ker V$ ($\tilde{M} = \mul\, V$, respectively). 
\item If $\dom\, V \subseteq X$ is closed and $\overline{\ran}\, V=Y$, then $V$ is a bounded operator. 
%\item If $\ran\, V \subseteq Y$ is closed and $\overline{\dom}\, V=X$, then $V^{-1}$ is a bounded operator.
\item If either $\left( \overline{\dom}\, V=X \wedge \ran \, V=Y\right)$ or $\left( \overline{\ran}\, V=Y \wedge \dom\, V=X\right) $ holds, then $V$ is a standard unitary operator.
\end{enumerate}
\end{proposition}
\textbf{Proof.}

(i) Let us first prove: $\tilde{M} = \mul\, V  \Rightarrow  M=\ker V $. 

Assume $\tilde{M} = \mul\, V $ and $ \left( {\begin{array}{*{20}c}
f\\
f'\\
\end{array} } \right)\in V$. Then
\[
 f\in M \Rightarrow \left( \left[ f,g \right]=\left[ f',g' 
\right]=0,\, \forall \left( {\begin{array}{*{20}c}
g\\
g'\\
\end{array} } \right)\in V \right) \Rightarrow \left( f'=0 \vee 0 \neq f'\in \mul\, V\right)  \Rightarrow f \in \ker V. 
\]
Hence, $M \subseteq \ker V$.
\[
f\in \ker V  \Rightarrow \left( \left[ f,g \right]=\left[ f',g' 
\right]=0,\, \forall \left( {\begin{array}{*{20}c}
g\\
g'\\
\end{array} } \right)\in V \right) \Rightarrow f \in M . 
\]
Hence, $\ker V \subseteq M$. This completes the proof of
\[
\tilde{M} = \mul\, V  \Rightarrow M=\ker V . 
\]
The converse $M = \ker V \Rightarrow \tilde{M}=\mul\, V$ follows immediately if we substitute the isometric $V^{-1}$ for the isometric $V$ in the previous implication. This completes the proof of (i). 

(ii) $V^{-1}\subseteq V^{\ast }$ and $\overline{\ran} \, V=Y$ imply $\overline{\dom}\, V^{\ast }=Y$. Then, according to Proposition \ref{proposition34} (iv) and (ii), we have $\mul\, \bar{V}=\left(\dom\, V^{\ast } \right)^{{[\bot ]}_{Y}}=\left\{ 0 \right\}$. Hence, $\bar{V}$ and $V$ are operators. 

Let us prove that $V$ is closed. Assume
\[
\left( {\begin{array}{*{20}c}
f_{n}\\
f_{n}'\\
\end{array} } \right)\in V\, \wedge \left( {\begin{array}{*{20}c}
f_{n}\\
f_{n}'\\
\end{array} } \right)\to \left( {\begin{array}{*{20}c}
f\\
f'\\
\end{array} } \right)\in \bar{V}\, \, (n\to \infty ).
\]
Because $\dom\, V$ is closed, $f\in \dom\, V$. Because $V$ is single-valued, it must be $\left( {\begin{array}{*{20}c}
f\\
f'\\
\end{array} } \right)\in V$. Hence, $V$ is closed. The continuity of $V$ follows from the closed graph theorem. 

(iii) Assume, for example, $\left( \overline{\dom}\, V=X \wedge \ran \, V=Y\right)$. This and $V^{-1}\subseteq V^{*}$ implies $V^{-1}= V^{*}$. Then statement (iv) follows from \cite[Proposition 2.3 (i)]{DHMS1}, or from \cite[Lemma 3.25]{W1}.
%Assume $\overline{\dom}\, V=X$ and $V$ is a surjective relation. Let us prove $V^{\ast }\subseteq V^{-1}$. 
%Let us arbitrarily select $\left( {\begin{array}{*{20}c}
%g'\\
%{0}\\
%\end{array} } \right)\in V^{\ast }\subseteq Y\times X$. Because $V^{-1}\subseteq V^{\ast }$, for $g'\in Y=\ran\, V$ there exists $g\in \dom\, V$ such that it holds $\left( {\begin{array}{*{20}c}
%g\\
%g'\\
%\end{array} } \right)\in V$. We need to prove $\left( {\begin{array}{*{20}c}
%g'\\
%{0}\\
%\end{array} } \right)\in V^{-1}$. Indeed, via definitions of isometric $V$ and adjoint $V^{\ast }$ relations, respectively, we get 
%\[
%\left[ f,g \right]_{X}=\left[ f',g' \right]_{Y}=\left[ f,g_{0} \right]_{X}, \forall \left( {\begin{array}{*{20}c}
%f\\
%f'\\
%\end{array} } \right)\in V.
%\]
%Because $\dom\, V$ is dense in $X$, we conclude $g_{0}=g$. Hence, $\left( 
%{\begin{array}{*{20}c}
%g'\\
%{0}\\
%\end{array} } \right)=\left( {\begin{array}{*{20}c}
%g'\\
%g\\
%\end{array} } \right)\in V^{-1}$. Together with $V^{-1}\subseteq V^{\ast }$ this gives $V^{-1}=V^{\ast }$. Hence, $V$ is a closed relation, and, according to Proposition \ref{proposition34} (v), $\ker V^{\ast }=\mul \, V=\left\{ 0 \right\}$, i.e. $V$ is an operator. 
\hfill $\square$
\\

Note that Proposition \ref{proposition42} (ii) is a generalization of \cite[Lemma 1.8.1]{BHS} because here we deal with arbitrary Krein spaces $X$ and $Y$ and we deal with a linear isometric relation $V$ with $\overline{\ran}\, V=Y$ rather than with a surjective operator $V$. 

Note also that in order to claim continuity of the operator $V$ when $X$ and $Y$ are Krein spaces, we had to assume $\overline{\ran}\, V=Y$ and $\dom\, V$ is closed, in Proposition \ref{proposition42} (ii). A more general statement holds when $X$ and $Y$ are Pontryagin spaces; see \cite[Lemma 2.1]{B3}. 

\section{Properties of the Green's boundary relation }\label{s10}

\textbf{5.1}. In Section \ref{s8}, we proved some general results about isometric relations between Krein spaces that we need to obtain the following properties of Green's boundary relations. 

\begin{proposition}\label{proposition52} Let $\Gamma : \mathcal{K}^{2}\to \mathcal{H}^{2}$ be a Green's boundary relation. Then:

\begin{enumerate}[(i)]%, (i), (ii),...
\item $\tilde{M} = \mul\, \Gamma \Leftrightarrow M=\ker \Gamma$.
%\item If $\ran\, \Gamma $ ($\dom\, \Gamma $) does not degenerate, then $M=\ker \Gamma $ ($\tilde{M} = \mul\, \Gamma $, respectively).
\item If $\dom\, \Gamma \subseteq \mathcal{K}^{2}$ is closed and $\overline{\ran}\, \Gamma=\mathcal{H}^{2}$, then $\Gamma$ is a bounded operator. 
%\item If $\ran\, \Gamma $ is closed and $\overline{\dom}\, \Gamma =\mathcal{K}^{2}$, then $\Gamma^{-1}$ is a bounded operator.
\item If either $\left( \overline{\dom}\,  \Gamma = \mathcal{K}^{2} \wedge \ran \,  \Gamma = \mathcal{H}^{2} \right)$ or $\left( \overline{\ran}\, \Gamma =\mathcal{K}^{2} \wedge \dom\,  \Gamma =\mathcal{K}^{2} \right) $ holds, then $ \Gamma $ is a standard unitary operator.
\end{enumerate}
\end{proposition}
\textbf{Proof.} Recall that the indefinite scalar products (\ref{eq28}) and (\ref{eq210}) define the Krein spaces $\left( \mathcal{K}^{2}, \left[ .,. \right]_{\mathcal{K}^{2}} \right)$, $(\mathcal{H}^{2}, \left[ .,. \right]_{\mathcal{H}^{2}})$, respectively. Since $\Gamma : \mathcal{K}^{2}\to \mathcal{H}^{2}$ satisfies (\ref{eq214}), $\Gamma$ is an isometric relation between those Krein spaces. We can substitute $\left( \mathcal{K}^{2},\, \left[ .,. \right]_{\mathcal{K}^{2}} \right)$, $(\mathcal{H}^{2},\, \left[ .,. \right]_{\mathcal{H}^{2}})$ and $\Gamma $ for $\left( X, {[.,.]}_{X} \right)$, $\left( Y,, {[.,.]}_{Y} \right)$ and $V$, respectively, in Proposition \ref{proposition42}. Also observe that the isotropic part of $\dom\, \Gamma $ in $\left( \mathcal{K}^{2},\, \left[ .,. \right]_{\mathcal{K}^{2}} \right)$ is $M=\dom\, \Gamma \cap \left(\dom\, \Gamma \right)^{+}$. Then statement (i)-(iii) are consequences of the corresponding statements in Proposition \ref{proposition42}.\hfill $\square$
\\

Note that Proposition \ref{proposition52} (ii) is a generalization of \cite[Lemma 1.8.1]{BHS} because here we allow $\mathcal{K}$ to be a Krein space, rather than a Hilbert space, and $\Gamma $ can be a linear relation with $\overline{\\ran}\, \Gamma= \mathcal{H}^{2}$ rather than a surjective operator.

\begin{corollary}\label{corollary56} Let $\Gamma :\mathcal{K}^{2}\to \mathcal{H}^{2}$ be a Green's boundary relation. If $\ran\, \Gamma$ is non-degenerate and condition (\ref{eq223}) is satisfied, then $M=N=S$.
\end{corollary}
\textbf{Proof.} According to Proposition \ref{proposition28} (ii), $M=S$. Since $\ran\, \Gamma$ is non-degenerate, $\tilde{M}=\lbrace0 \rbrace$. Because $\mul\, \Gamma \subseteq \tilde{M}$, $\Gamma$ is an operator and it holds $\tilde{M}=\mul\, \Gamma $ . This and Proposition \ref{proposition52} (i) imply $\ker \, \Gamma = M $. Hence, $ M=N=S$. \hfill $\square$ 

\begin{corollary}\label{corollary58} Let $\Gamma :\mathcal{K}^{2} \to \mathcal{H}^{2}$ be a Green's boundary operator. Then, $\dom\, \Gamma$ is closed, $\ran\, \Gamma=\mathcal{H}^{2}$ and (\ref{eq222}) holds if and only if $\Pi=(\mathcal{H}, \Gamma_{0}, \Gamma_{1})$ is an ordinary boundary triple for $S^{+}:=\dom\, \Gamma $.
\end{corollary}
\textbf{Proof.} $\dom\, \Gamma$ is closed and (\ref{eq222}) holds imply that (\ref{eq223}) holds. This means that $S:=\left( \dom\, \Gamma\right)^{+} $ is symmetric. According to Proposition \ref{proposition52} (ii), $\Gamma$ is a continuous operator. Because, it is also surjective, the triple $\Pi=(\mathcal{H}, \Gamma_{0}, \Gamma_{1})$ is ordinary.

Conversely, if $\Pi=(\mathcal{H}, \Gamma_{0}, \Gamma_{1})$ is an ordinary boundary triple for $S^{+}$, then $\dom\, \Gamma = S^{+}$ is closed, $\ran\, \Gamma =\mathcal{H}^{2}$. The existence of closed symmetric relation $S$ guaranties that (\ref{eq222}) holds. \hfill $\square$ 
\\

Similar to the notation $M,\, N, \, S$ related to $\dom\, \Gamma $, we introduce the notation $\tilde{M}:= \ran\, \Gamma \cap \left( \ran\, \Gamma \right)^{[\bot]}$, $\tilde{N}:= \mul\, \Gamma$, and $\tilde{S}:= \left( \ran\, \Gamma\right)^{+}$ related to $\ran\, \Gamma$. 
\begin{proposition}\label{proposition510}. Let $\Gamma : \mathcal{K}^{2} \rightarrow  \mathcal{H}^{2}$ be a 
Green's boundary relation. Then:
\begin{enumerate}[(i)]%, (i), (ii),...
\item $\Gamma $ is a unitary relation if and only if $\Gamma $ satisfies the condition
\begin{equation}
\label{eq54}
\left( \left[ \hat{f},\hat{g} \right]_{\mathcal{K}^{2}}=\left[ \hat{h},\hat{k} 
\right]_{\mathcal{H}^{2}}, \forall \left( {\begin{array}{*{20}c}
\hat{f}\\
\hat{h}\\
\end{array} } \right)\in \Gamma \right)\Rightarrow \left( 
{\begin{array}{*{20}c}
\hat{g}\\
\hat{k}\\
\end{array} } \right)\in \Gamma .
\end{equation}
\item (\ref{eq54}) $\Rightarrow $ (\ref{eq223}).
\item If $\Gamma$ is a unitary boundary relation, then $M=N=S$ and  $\tilde{M}=\tilde{N}=\tilde{S}$.
\end{enumerate}
\end{proposition}
\textbf{Proof.} 

(i) We need to prove that for a GBR $\Gamma$, it holds $(\ref{eq54}) \Leftrightarrow \Gamma^{-1}=\Gamma^{*} $. 

If we substitute $\mathcal{K}^{2}$ for $X$, $\mathcal{H}^{2}$ for $Y$, and $\Gamma$ for $V$, then $ \left(  X\times Y, \left[ .,. \right]_{X \times Y} \right) = \left(  \mathcal{K}^{2}\times \mathcal{H}^{2}, \left[ .,. \right]_{\mathcal{K}^{2}\times \mathcal{H}^{2}} \right) $. Therefore, $ (\ref{eq54}) \Leftrightarrow\Gamma ^{\left[ \perp \right]_{\mathcal{K}^{2}\times \mathcal{H}^{2}}} \subseteq \Gamma $. 

Because, $\Gamma$ is a GBR, hence isometric, according to Proposition \ref{proposition38}, it holds $\Gamma  \subseteq \Gamma ^{\left[ \perp \right]_{\mathcal{K}^{2}\times \mathcal{H}^{2}}} $.  

Therefore, for a GBR $\Gamma$, $(\ref{eq54}) \Leftrightarrow \Gamma ^{\left[ \perp \right]_{\mathcal{K}^{2}\times \mathcal{H}^{2}}} = \Gamma $. According to the unitary part of Proposition \ref{proposition38}, this is equivalent to $\Gamma^{-1}=\Gamma^{*}$. This proves (i).

(ii) Assume that (\ref{eq54}) holds.
\[
\left( \left[ f',g \right]-\left[ f,g' \right]=0, \forall \hat{f}\in 
\dom\, \Gamma \right)\Rightarrow \left( \left[ \hat{f},\hat{g} 
\right]_{\mathcal{K}^{2}}=\left[ \hat{h},\hat{k} \right]_{\mathcal{H}^{2}}=0, \forall \left( 
{\begin{array}{*{20}c}
\hat{f}\\
\hat{h}\\
\end{array} } \right)\in \Gamma \right).
\]
According to (\ref{eq54}), $\left( {\begin{array}{*{20}c}
\hat{g}\\
\hat{k}\\
\end{array} } \right)\in \Gamma $. Hence, $\hat{g}\in \dom\, \Gamma$. This means that (\ref{eq223}) holds.

(iii) If $\Gamma$ is a unitary relation, then according to (i) and (ii), $\Gamma$ satisfies (\ref{eq223}). This implies $S=M$. 

If we substitute $\Gamma$ for $H$ in Proposition \ref{proposition34} (v) and use $\Gamma^{-1}=\Gamma^{*}$, we get: $\ker\, \Gamma^{*} =\mul \, \Gamma = (\ran\, \Gamma)^{[\bot]_{{\mathcal{H}}^{2}}}$. This and $ \mul\, \Gamma \subseteq \dom\, \Gamma^{*} = \ran\, \Gamma$ implies $\mul\, \Gamma =\tilde{M}$. According to Proposition \ref{proposition52} (i), it follows $M=\ker \, \Gamma=N$. Therefore, $M=N=S$.

If we substitute $\Gamma^{-1}$ for $\Gamma$ in previous reasoning we conclude $\tilde{M}=\tilde{N}=\tilde{S}$.\hfill $\square$
\\

\textbf{5.2.} In this subsection, we will examine some consequences of the main transformation; refer to \cite[Proposition 2.10]{DHMS1} and \cite[(3.5)]{BDHS}. Let us consider the spaces $(\mathcal{K}^{2},\left[ .,. \right]_{\mathcal{K}^{2}})$, $(\mathcal{H}^{2},\left[ .,. \right]_{\mathcal{H}^{2}})$ and $\mathcal{K}\times \mathcal{H}$. Elements of $\mathcal{K}^{2}$ are denoted by $\hat{f}=\left( {\begin{array}{*{20}c}
f\\
f'\\
\end{array} } \right),\, \hat{g}=\left( {\begin{array}{*{20}c}
g\\
g'\\
\end{array} } \right)$, and elements of $\mathcal{H}^{2}$ by $\hat{h}=\left( 
{\begin{array}{*{20}c}
h\\
h'\\
\end{array} } \right), \hat{k}=\left( {\begin{array}{*{20}c}
k\\
k'\\
\end{array} } \right)$. Then elements of $\mathcal{K}\times \mathcal{H}$ are $\left( 
{\begin{array}{*{20}c}
f\\
h\\
\end{array} } \right), \left( {\begin{array}{*{20}c}
f'\\
h'\\
\end{array} } \right), \left( {\begin{array}{*{20}c}
g\\
k\\
\end{array} } \right), \left( {\begin{array}{*{20}c}
g'\\
k'\\
\end{array} } \right)$ etc. Let the linear mapping $\mathcal{J}:\mathcal{K}^{2}\times \mathcal{H}^{2}\to \left( \mathcal{K}\times \mathcal{H} \right)\times \left( \mathcal{K}\times \mathcal{H} \right)$ be the well known \textit{main transformation}; see e.g. \cite[(2.16)]{DHMS1}. It satisfies 
\begin{equation}
\label{eq56}
\mathcal{J}:\Gamma = \left\{\left( {\begin{array}{*{20}c}
f\\
f'\\
\end{array} } \right ), \left( {\begin{array}{*{20}c}
h\\
h'\\
\end{array} } \right)\right\} \rightarrow\tilde{A}=\left\{\left( {\begin{array}{*{20}c}
f\\
h\\
\end{array} } \right ), \left( {\begin{array}{*{20}c}
f'\\
-h'\\
\end{array} } \right)\right\}  
\end{equation}
Recall that the product Krein space $\mathcal{K}\times \mathcal{H}$ is endowed with the inner product 
\begin{equation}
\label{eq58}
\left[ \left( {\begin{array}{*{20}c}
f\\
h\\
\end{array} } \right),\left( {\begin{array}{*{20}c}
g\\
k\\
\end{array} } \right) \right]_{\mathcal{K}\times \mathcal{H}}:=\left[ f,g \right]_{\mathcal{K}}+\left( h,k 
\right)_{\mathcal{H}}.
\end{equation}
and that $\tilde{A}:=\mathcal{J}(\Gamma)$ is a symmetric (self-adjoint) relation if and only if $U$ is isometric (unitary) relation.

%Note that the main transformation $\mathcal{J}$ is a unitary operator, i.e. it satisfies $\mathcal{J}^{*}=\mathcal{J}^{-1}$. 

The assumption $\mathcal{K}_{1}\cap \mathcal{K}_{2}=\lbrace 0 \rbrace$ in the following proposition means that $\mathcal{K}_{1}$ and $\mathcal{K}_{2}$ are distinct vector spaces; for example, they cannot be the same vector spaces with different scalar products. We will present the proof of the following lemma because we will need it in the examples that follow to identify the reasons why the proof does not apply to those examples. 

\begin{lemma}\label{lemma518} \cite[Lemma 2.1]{B1} Assume that $\mathcal{K}_{1}$ and $\mathcal{K}_{2}$ are Krein spaces, $\mathcal{K}_{1}\cap \mathcal{K}_{2}=\lbrace 0 \rbrace$, and $A_{i}\subseteq \mathcal{K}_{i}\times\mathcal{K}_{i},\thinspace i=1,\thinspace 2$, are linear relations. We can define direct orthogonal sum
\[
\tilde{\mathcal{K}}:=\mathcal{K}_{1}\left[ + \right]\mathcal{K}_{2}, 
\]
and 
\begin{equation}
\label{eq518}
\tilde{A}:=A_{1}\left[ + \right]A_{2}:=\left\{ \left\{\left( {\begin{array}{*{20}c}
f_{1}\\
f_{2}\\
\end{array} } \right ), \left( {\begin{array}{*{20}c}
f'_{1}\\
f'_{2}\\
\end{array} } \right)\right\} : \left( {\begin{array}{*{20}c}
f_{i}\\
f_{i}'\\
\end{array} } \right)\in A_{i},\thinspace  i=1,\thinspace 2 \right\}\subseteq \tilde{\mathcal{K}}\times \tilde{\mathcal{K}}.
\end{equation}
%where $A_{i} \subseteq \tilde{A},\, i=1,\, 2$.  

Then linear relation $\tilde{A}:=A_{1}\left[ + \right]A_{2}$ is symmetric (self-adjoint) in $\tilde{\mathcal{K}}$ if and only if the linear relations $A_{i}\subseteq \mathcal{K}_{i}\times \mathcal{K}_{i}, \thinspace i=1,\thinspace 2$, are symmetric (self-adjoint). 
\end{lemma} 
\textbf{Proof.} First, assume that $A_{i},\,  i=\,1,\, 2$ are symmetric and verify $A_{1}\left[ + \right]A_{2} \subseteq \left( A_{1}\left[ + \right]A_{2}\right)^{+}$. It holds
\[
\left\{\left( {\begin{array}{*{20}c}
g_{1}[+]g_{2}\\
g_{1}'[+]g_{2}'\\
\end{array} } \right )\right\}  =
\left\{\left( {\begin{array}{*{20}c}
g_{1}\\
g_{2}\\
\end{array} } \right ), \left( {\begin{array}{*{20}c}
g'_{1}\\
g'_{2}\\
\end{array} } \right)\right\} \in \tilde{A}
\]
\[
\Rightarrow \left[ \left( {\begin{array}{*{20}c}
f_{1}\\
f_{2}\\
\end{array} } \right ), \left( {\begin{array}{*{20}c}
g_{1}'\\
g_{2}'\\
\end{array} } \right)\right] - \left[ \left( {\begin{array}{*{20}c}
f'_{1}\\
f'_{2}\\
\end{array} } \right ), \left( {\begin{array}{*{20}c}
g_{1}\\
g_{2}\\
\end{array} } \right)\right]=0, \forall \left\{\left( {\begin{array}{*{20}c}
f_{1}\\
f_{2}\\
\end{array} } \right ), \left( {\begin{array}{*{20}c}
f'_{1}\\
f'_{2}\\
\end{array} } \right)\right\} \in \tilde{A}
\]
\[
\Rightarrow \left\{\left( {\begin{array}{*{20}c}
g_{1}\\
g_{2}\\
\end{array} } \right ), \left( {\begin{array}{*{20}c}
g'_{1}\\
g'_{2}\\
\end{array} } \right)\right\} \in \tilde{A}^{+}. 
\]

Conversely, assume $A_{1}\left[ + \right]A_{2} \subseteq \left( A_{1}\left[ + \right]A_{2}\right)^{+}$. It holds
\[
A_{i} \subseteq  A_{1}\left[ + \right]A_{2} \subseteq \left( A_{1}\left[ + \right]A_{2}\right)^{+} \subseteq A_{i}^{+},\, i=\, 1,\, 2.
\]

Now assume $\left( A_{1}\left[ + \right]A_{2}\right)^{+} = A_{1}\left[ + \right]A_{2} $ and prove that both $A_{i}, \, i=1,2$, are self-adjoint. We have 
\[ 
\left( {\begin{array}{*{20}c}
k_{1} \\
k_{1}'\\
\end{array} } \right ) \in A_{1}^{+} \Rightarrow \left[ h_{1},k_{1}' \right] = \left[ h_{1}',k_{1} \right], 
\forall \left( {\begin{array}{*{20}c}
h_{1} \\
h_{1}'\\
\end{array} } \right ) \in A_{1} 
\]
\[ 
\Rightarrow  \left[ h_{1}[+]h_{2},k_{1}' \right] = \left[ h'_{1}[+]h'_{2},k_{1} \right], 
\forall \left( {\begin{array}{*{20}c}
h_{1}[+]h_{2} \\
h_{1}'[+]h_{2}'\\
\end{array} } \right ) \in A_{1}\left[ + \right]A_{2}
\]
\[
\Rightarrow \left( {\begin{array}{*{20}c}
k_{1} \\
k_{1}'\\
\end{array} } \right ) \in \left( A_{1}\left[ + \right]A_{2}\right)^{+}=A_{1}\left[ + \right]A_{2} \Rightarrow 
\left( {\begin{array}{*{20}c}
k_{1} \\
k_{1}'\\
\end{array} } \right ) \in A_{1}. 
\]
This proves $ A_{1}^{+} \subseteq  A_{1}$. By the same token it holds $A_{2}^{+} = A_{2}$. 

Conversely, we assume that $A_{1}$ and $A_{2}$ are self-adjoint. Then
\[
\left\{\left( {\begin{array}{*{20}c}
g_{1}\\
g_{2}\\
\end{array} } \right ), \left( {\begin{array}{*{20}c}
g'_{1}\\
g'_{2}\\
\end{array} } \right)\right\} \in \tilde{A}^{+} 
\]
\[
\Rightarrow \left[ \left( {\begin{array}{*{20}c}
f_{1}\\
f_{2}\\
\end{array} } \right ), \left( {\begin{array}{*{20}c}
g_{1}'\\
g_{2}'\\
\end{array} } \right)\right] - \left[ \left( {\begin{array}{*{20}c}
f'_{1}\\
f'_{2}\\
\end{array} } \right ), \left( {\begin{array}{*{20}c}
g_{1}\\
g_{2}\\
\end{array} } \right)\right]=0, \forall \left\{\left( {\begin{array}{*{20}c}
f_{1}\\
f_{2}\\
\end{array} } \right ), \left( {\begin{array}{*{20}c}
f'_{1}\\
f'_{2}\\
\end{array} } \right)\right\} \in \tilde{A}.
\]
If we select $f_{2}=f'_{2}=0$, it follows 
\[
 \left[ f_{1},g'_{1} \right] - \left[ f'_{1},g_{1} \right]=0,\,   \forall \left( {\begin{array}{*{20}c}
f_{1}\\
f'_{1}\\
\end{array} } \right ) \in A_{1}.
\]
Because,  we assume $A_{1}$ is self-adjoint, it follows $\left( {\begin{array}{*{20}c}
g_{1}\\
g'_{1}\\
\end{array} } \right ) \in A_{1}$. 
In the same way, $\left( {\begin{array}{*{20}c}
g_{2}\\
g'_{2}\\
\end{array} } \right ) \in A_{2}$. Therefore, $\left\{\left( {\begin{array}{*{20}c}
g_{1}\\
g_{2}\\
\end{array} } \right ), \left( {\begin{array}{*{20}c}
g'_{1}\\
g'_{2}\\
\end{array} } \right)\right\} \in A_{1}[+]A_{2}$.

This proves $\tilde{A}^{+}=\tilde{A}$. \hfill $\square$ 
\\

The assumption $\mathcal{K}_{1}\cap \mathcal{K}_{2}=\lbrace 0 \rbrace $ is \textbf{not} satisfied, for example, when $\mathcal{K}_{1}= \mathcal{K}_{2}=\mathbb{C}$, but $\mathcal{K}_{1}$ and $\mathcal{K}_{2}$ are different Krein spaces because they are endowed with different scalar products. This assumption is essential because without it, we would not be able to prove the implication: $\left( A_{1}\left[ + \right]A_{2}\right)^{+}=A_{1}\left[ + \right]A_{2} \Rightarrow  A_{i}^{+} = A_{i}, \, i=1.2$. Namely, if addition  $h_{1}+h_{2} $ were defined in the previous proof, we would not be able to conclude $A_{i}^{+} \subseteq \mathcal{K}_{i}, \, i=1,2$. 

Recall that the linear mapping $\mathcal{J}^{-1}$ maps a self-adjoint relation $\tilde{A}$ in $ \mathcal{K} [+] \mathcal{H}$ to a unitary relation $\Gamma :\mathcal{K}^{2}\rightarrow \mathcal{H}^{2}$, i.e. 
\begin{equation}
\label{eq522}
\mathcal{J}^{-1}: \tilde{A}=\left\{\left( {\begin{array}{*{20}c}
f\\
h\\
\end{array} } \right ), \left( {\begin{array}{*{20}c}
f'\\
h'\\
\end{array} } \right)\right\} \rightarrow \Gamma = \left\{\left( {\begin{array}{*{20}c}
f\\
f'\\
\end{array} } \right ), \left( {\begin{array}{*{20}c}
h\\
-h'\\
\end{array} } \right)\right\}.
\end{equation}
By definition; see \cite[Defintion 2.2]{B1}, the linear relation $\tilde{A}$ that satisfies conditions of Lemma \ref{lemma518} is reducible. 

Now we have everything needed to prove that the self-adjoint linear relation $\tilde{A}$ is reduced by the decomposition $ \tilde{\mathcal{K}}=\mathcal{K} [+] \mathcal{H}$ if and only if $\Gamma=\mathcal{J}^{-1}\left( \tilde{A}\right)$ is a trivial unitary relation. 

\begin{proposition}\label{proposition522} Let $ \tilde{A}$ be a self-adjoint linear relation that satisfies the conditions of Lemma \ref{lemma518} with $\mathcal{K}=\mathcal{K}_{1}$ and $\mathcal{H}=\mathcal{K}_{2}$. Then the relation $\Gamma=\mathcal{J}^{-1}\left( \tilde{A}\right) : \mathcal{K}^{2}\rightarrow \mathcal{H}^{2}$ is a \textbf{trivial} unitary Green's boundary relation. 

Conversely, let $\Gamma: \mathcal{K}^{2}\rightarrow \mathcal{H}^{2}$ be a \textbf{trivial} unitary Green's boundary relation. Then $\tilde{A}=\mathcal{J}\left( \Gamma\right)$ is a self-adjoint relation that satisfies the conditions of Lemma \ref{lemma518} with $\mathcal{K}=\mathcal{K}_{1}$ and $\mathcal{H}=\mathcal{K}_{2}$.
\end{proposition}
\textbf{Proof.}  According to the condition $A_{i}\subseteq \mathcal{K}_{i}\times\mathcal{K}_{i},\thinspace i=1,\thinspace 2$, and according (\ref{eq522}), $A_{1} = \dom\, \Gamma $ and $ A_{2} = \ran \, \Gamma $. Because the self-adjoint relation $ \tilde{A}$ satisfies the conditions of Lemma \ref{lemma518}, $A_{1} = \dom\, \Gamma $ and $ A_{2} =\ran \, \Gamma $ are self-adjoint relations. Then $\dom\, \Gamma \subseteq \mathcal{K}^{2}$ and $\ran \, \Gamma $ are hyper-maximal neutral i.e., $\dom\, \Gamma = M$ and $\ran \, \Gamma = \tilde{M}$. 

Then according to Proposition \ref{proposition510} (iii), it follows $\dom\, \Gamma = \ker\, \Gamma $ and $\ran \, \Gamma = \mul\, \Gamma $. 

Conversely, if $\Gamma$ is a trivial unitary relation, then $\tilde{A}=\mathcal{J}\left( \Gamma\right)$ is a self-adjoint relation, and $\dom\, \Gamma = \ker\, \Gamma $ and $\ran \, \Gamma = \mul\, \Gamma $. According to Proposition \ref{proposition510} (iii), it follows $\dom\, \Gamma = M $ and $\ran \, \Gamma = \tilde{M} $.  This further means that $A_{1}:=\dom\, \Gamma \subseteq \mathcal{K}^{2} $ and $A_{2}:=\ran \, \Gamma \subseteq \mathcal{H}^{2}$ are self-adjoint relations in their respective Krein spaces, and $\tilde{A}:=A_{1}[+]A_{2}$. \hfill $\square$ 
\\

At the end of the section, we provide two examples of the irreducible symmetric linear relations in the Krein space of the form $\tilde{\mathcal{K}}=\mathcal{K}[+]\mathcal{H}$. In the first example, we demonstrate that even when $\mathcal{K}_{1}\cap \mathcal{K}_{2}=\lbrace 0 \rbrace$, the symmetric relation $\tilde{A}$ may not have the decomposition (\ref{eq518}) with symmetric relations $A_{i},\, i=1,2$. 

\begin{example}\label{example520} Let $\mathcal{K}$, $\mathcal{H}$ and $\Gamma$ be defined by (\ref{eq14}) and (\ref{eq16}). 
\end{example}

Then $\mathcal{K}\cap \mathcal{H}= L^{2}(0,\infty ) \cap \mathbb{C}=\lbrace 0 \rbrace$ and $\Gamma: \mathcal{K}^{2}\rightarrow \mathcal{H}^{2}$ is isometric. We know that $\Pi=(\mathcal{H}, \Gamma_{0}, \Gamma_{1})$ is an ordinary boundary triple.  We also know that the main transformation of $\Gamma$, $\tilde{A}:= \mathcal{J}\left( \Gamma \right): \mathcal{K}\times  \mathcal{H} \rightarrow \mathcal{K}\times \mathcal{H} $, is a symmetric relation. According to (\ref{eq56}), it holds
\[
\tilde{A}=\mathcal{J}(\Gamma)=\left\{\left\{\left( {\begin{array}{*{20}c}
f\\
\Gamma_{0} \hat{f}\\
\end{array} } \right ), \left( {\begin{array}{*{20}c}
f'\\
-\Gamma_{1} \hat{f}\\
\end{array} } \right)\right\} : \left( {\begin{array}{*{20}c}
f\\
f'\\
\end{array} } \right ) \in \dom\, \Gamma\right\}.
\] 
Assume that the decomposition (\ref{eq518}), $\tilde{A}=A_{1}[+]A_{2}$, exists. Then according to the assumption  $A_{i}\subseteq \mathcal{K}_{i}\times\mathcal{K}_{i},\thinspace i=1,\thinspace 2$, it would be
\[
\tilde{A} \cap \left( \mathcal{K}[+]\lbrace 0 \rbrace \right) =\left\{\left\{\left( {\begin{array}{*{20}c}
f\\
0\\
\end{array} } \right ), \left( {\begin{array}{*{20}c}
f'\\
0\\
\end{array} } \right)\right\} : \left( {\begin{array}{*{20}c}
f\\
f'\\
\end{array} } \right ) \in \ker \Gamma\right\} \Rightarrow  A_{1}= \ker \Gamma \neq \dom\, \Gamma.
\]
Similarly $\tilde{A} \cap \left( \lbrace 0 \rbrace [+] \mathcal{H}\right) = \lbrace 0 \rbrace = A_{2}$. 
This means that the symmetric linear relation $\tilde{A}$ in this example cannot be reduced by symmetric components in the sense of Lemma \ref{lemma518} because the condition $A_{i}\subseteq \mathcal{K}_{i}\times\mathcal{K}_{i},\thinspace i=1,\thinspace 2$, is not satisfied. \hfill $\square$

\begin{example}\label{example519} Let the vector spaces $\mathcal{K}_{1}$ and $\mathcal{K}_{2}$, $\mathcal{K}_{1}=\mathcal{K}_{1}=\mathbb{C}$, be endowed with (different) scalar products $\left[.,. \right]$ and $\left(.,. \right)$, respectively, that satisfy
\[
\left[f,g \right] =-\left(f,g \right) :=-f\overline{g}, \, \forall f,g \in \mathbb{C}. 
\]
We can consider that $\mathcal{K}_{1}$ and $\mathcal{K}_{2}$ are Krein and Hilbert space. Then the linear relations
\[
A_{1}=A_{2}=A := \left\{\left( {\begin{array}{*{20}c}
f\\
f'\\
\end{array} } \right ) :\,   f, f'\in  \mathbb{C}\right\},
\]
satisfy the condition $A_{i}\subseteq \mathcal{K}_{i} \times \mathcal{K}_{i}, \,i=1, 2$, but they are not self-adjoint because $A_{i}^{+}=\lbrace 0 \rbrace,\, i=1,\, 2$. Let us prove that the relation $\tilde{A}$ defined by
\begin{equation}
\label{eq520}
\tilde{A}= A \left[ + \right]A:=\left\{\left( {\begin{array}{*{20}c}
f\\
f\\
\end{array} } \right ), \left( {\begin{array}{*{20}c}
f'\\
f'\\
\end{array} } \right)\right\},
\end{equation}
is a self-adjoint relation in $\mathcal{K}_{1} [+] \mathcal{K}_{2}$ endowed with the scalar product (\ref{eq58}).
\end{example} 

Indeed, for 
\[
\left\{\left( {\begin{array}{*{20}c}
g_{1}\\
g_{2}\\
\end{array} } \right ), \left( {\begin{array}{*{20}c}
g_{1}'\\
g_{2}'\\
\end{array} } \right)\right\} \in \left( \tilde{A}\right)^{+} 
\]
\[
\Rightarrow \left[ \left( {\begin{array}{*{20}c}
f\\
f\\
\end{array} } \right ), \left( {\begin{array}{*{20}c}
g_{1}'\\
g_{2}'\\
\end{array} } \right)\right] - \left[ \left( {\begin{array}{*{20}c}
f'\\
f'\\
\end{array} } \right ), \left( {\begin{array}{*{20}c}
g_{1}\\
g_{2}\\
\end{array} } \right)\right]=
\]
\[
=-f\overline{g_{1}'}+ f\overline{g_{2}'}-\left(-f'\overline{g_{1}}+ f'\overline{g_{2}}\right)=0, \forall \left\{\left( {\begin{array}{*{20}c}
f\\
f\\
\end{array} } \right ), \left( {\begin{array}{*{20}c}
f'\\
f'\\
\end{array} } \right)\right\} \in \tilde{A}.
\]
Because, $f \in \mathbb{C}$ and $f' \in \mathbb{C}$ are by assumption not identically zero, it follows 
\[
 g_{1}=g_{2}=:g \, \wedge \, g'_{1}=g'_{2}=:g'
\]
\[
\Rightarrow A^{+}=\left\{\left( {\begin{array}{*{20}c}
g\\
g\\
\end{array} } \right ), \left( {\begin{array}{*{20}c}
g'\\
g'\\
\end{array} } \right)\right\} = \tilde{A}.
\]

Hence, $\tilde{A}$ is self-adjoint.  \hfill $\square$
\section{An application of the Green's boundary relations on Weyl families}\label{s12}
In this section, we use the concept of the Green's boundary relation to prove generalizations of \cite[Lemma 3.15]{BDHS} and \cite[Theorem 4.8]{BDHS}.  

For a linear relation $T\subseteq \mathcal{K}^{2}$, the following concepts are defined:
\[
R_{z}\left( T \right):=\ker \left( T-z \right), \hat{R}_{z}\left( T 
\right):=\left\{ \left( {\begin{array}{*{20}c}
f_{z}\\
zf_{z}\\
\end{array} } \right):f_{z}\in R_{z}\left( T \right) \right\}, z\in \mathbb{C}.
\]
Let $\Gamma :\mathcal{K}^{2} \to \mathcal{H}^{2}$ be a Green's boundary relation. If we denote 
$T:=\dom\, \Gamma $, then for $\hat{f_{z}}:= \left( {\begin{array}{*{20}c}
f_{z}\\
zf_{z}\\
\end{array} } \right) \in \hat{R}_{z}\left( T 
\right)$ and $\hat{h}=\left( {\begin{array}{*{20}c}
h\\
h'\\
\end{array} } \right)\in \mathcal{H}^{2}$, the relation 
%\begin{equation}
%\label{eq62}
%\gamma \left( z \right)=\left\{ \left( {\begin{array}{*{20}c}
%h\\
%f_{z}\\
%\end{array} } \right):\left( {\begin{array}{*{20}c}
%\hat{f_{z}}\\
%\hat{h}\\
%\end{array} } \right)\in \Gamma ,\, \hat{f_{z}}\in \hat{R}_{z}\left( T \right) \right\}
%\end{equation}
%then the $\gamma $\textit{-field and the Weyl family associated with the boundary relation} $\Gamma $ are defined %analogously as $\gamma $-field and the Weyl function in the case of ordinary boundary triple, for example the relation 
\begin{equation}
\label{eq64}
M\left( z \right)=\left\{ \hat{h}:\left( {\begin{array}{*{20}c}
\hat{f_{z}}\\
\hat{h}\\
\end{array} } \right)\in \Gamma :\hat{f_{z}}\in \hat{R}_{z}\left( T \right) 
\right\}=\Gamma \left( \hat{R}_{z}\left( T \right) \right),
\end{equation}
is termed the \textit{Weyl family associated with the boundary relation} $\Gamma $ (cf. \cite[Definition 3.6]{BDHS}).

The following definition is a generalization of \cite[Definition 3.11]{BDHS} and \cite[Definition 3.4]{DHMS1}. In those definitions, $\Gamma$ is assumed to be a unitary relation, while in the following definition we assume only that the weaker condition (\ref{eq222}) is satisfied. 

\begin{definition}\label{definition62} A Green's boundary relation $\Gamma :\mathcal{K}^{2}\to \mathcal{H}^{2}$ that satisfies condition (\ref{eq222}), with $T:=\dom\, \Gamma $, is called minimal if 
\begin{equation}
\label{eq66}
\mathcal{K}:=clos\left\{ R_{z}\left( T \right):z\in \hat{\rho }\left( S \right) \right\},
\end{equation}
where $S:=\left( \dom\, T\right)^{+}$. 
\end{definition} 
It is necessary to have $\hat{\rho }\left( S \right)\ne \emptyset $; otherwise it would be $\mathcal{K}=\emptyset$.

The following proposition is a generalization of \cite[Lemma 3.15]{BDHS} because here we do not assume that $\Gamma$ is a unitary linear relation, as was implicitly assumed in \cite[Definition 3.11]{BDHS}. Here we assume only that the weaker condition (\ref{eq222}) holds. The proof is the same as the proof of Lemma \cite[Lemma 3.15]{BDHS}. We will repeat the proof here to demonstrate that it holds in this more general situation.  
\begin{proposition}\label{proposition64} Let $\mathcal{K}$ be a Pontryagin space and let $\Gamma :\mathcal{K}^{2}\to \mathcal{H}^{2}$ be a minimal Green's boundary relation. Then $S$ is a symmetric operator with $\sigma_{p}(S)=\emptyset$.
\end{proposition}
\textbf{Proof}. Because the minimal Green's boundary relation $\Gamma$ satisfies condition (\ref{eq222}), according to Proposition \ref{proposition28} (i), $S$ is a closed symmetric linear relation in the Pontryagin space $\mathcal{K}$. The following is a proof that the relation $S$ does not have any eigenvalues in $\mathbb{C}\cup\lbrace\infty\rbrace$. 

Assume that the claim is not true. Let us first assume that there exist a finite $\lambda \in \sigma_{p}(S)$ with $\left( {\begin{array}{*{20}c}
f_{\lambda}\\
\lambda f_{\lambda}\\
\end{array} } \right) \in S$.  Let $\eta \in \hat{\rho }\left( S \right)$ be an arbitrarily selected complex number and $g_{\eta} \in R_{\eta}(T):=\ker \left( T-\eta I \right) $, i.e. $\left( {\begin{array}{*{20}c}
g_{\eta}\\
\eta g_{\eta}\\
\end{array} } \right) \in T  $. Because $S=T^{+}$ it follows
\[\left[ \left( {\begin{array}{*{20}c}
f_{\lambda}\\
\lambda f_{\lambda}\\
\end{array} } \right), \left( {\begin{array}{*{20}c}
g_{\eta}\\
\eta g_{\eta}\\
\end{array} } \right) \right]_{\mathcal{K}^{2}}=0 \Rightarrow \left( \bar{\eta}-\lambda\right)\left[f_{\lambda},g_{\eta} \right] =0 \Rightarrow f_{\lambda} \left[ \perp \right]_{\mathcal{K}} R_{\eta}(T).
\]
Because $\eta \in \hat{\rho }\left( S \right)$ was arbitrarily selected, according to (\ref{eq66}), $f_{\lambda}=0$, i.e., $\lambda$ is not an eigenvalue of $S$. 

The assumption $\lambda = \infty$ is an eigenvalue of $S$ means $\left( {\begin{array}{*{20}c}
0\\
f_{\infty}\\
\end{array} } \right) \in S$. We have 
\[\left[ \left( {\begin{array}{*{20}c}
0\\
f_{\infty}\\
\end{array} } \right), \left( {\begin{array}{*{20}c}
g_{\eta}\\
\eta g_{\eta}\\
\end{array} } \right) \right]_{\mathcal{K}^{2}}=0 \Rightarrow \left[f_{\infty},g_{\eta} \right] =0
\Rightarrow f_{\infty} \left[ \perp \right]_{\mathcal{K}} R_{\eta}(T).
\] 
Again, according to (\ref{eq66}), it holds $f_{\infty}=0$, which means that $S$ is an operator. 
\hfill $\square$
\\

The following theorem is a generalization of \cite[Theorem 4.8]{BDHS} because we do not assume the existence of the closed symmetric linear relation $S$ as was assumed in \cite[Theorem 4.8]{BDHS}.

\begin{theorem}\label{theorem66} Let $\mathcal{K}$ be a Pontryagin space with $\kappa $ negative squares. Let $\Gamma :\mathcal{K}^{2}\to \mathcal{H}^{2}$ be a minimal unitary boundary relation, and let $M\left( z \right)$ be the Weyl family associated with $\Gamma $, i.e., the Weyl family defined by (\ref{eq64}). Then the family $M\left( z \right)$ is a generalized Nevanlinna family with $\kappa $ negative squares. 
\end{theorem}
\textbf{Proof.} According to Proposition \ref{proposition510}, the unitary relation $\Gamma$ satisfies condition (\ref{eq222}). Then, according to Proposition \ref{proposition28} (i), $S:=\left( \dom\, \Gamma\right)^{+}=T^{+}$ is a closed symmetric linear relation in the Pontryagin space $\mathcal{K}$ that satisfies $S^{+}=\overline{\dom}\, \Gamma$. This means that the \textbf{unitary boundary relation} $\Gamma $ in the sense of Definition \ref{definition214} is a \textbf{boundary relation} $\Gamma $ for $S^{+}$ in terms of \cite{BDHS}. Then Theorem \ref{theorem66} follows from \cite[Theorem 4.8]{BDHS}.  \hfill $\square$ 

\section{Green's boundary relations with a non-degenerate range }\label{s14}

In this section we investigate Green's boundary relations such that $\ran\, \Gamma$ is a non-degenerate manifold. Sometimes, we use the stronger assumption $\ran\, \Gamma = \mathcal{H}^{2}$. In the latter case, the relation $\Gamma$ becomes a \textit{Green's boundary operator}. As usual, $\bar{\Gamma}$ denotes the closure of $\Gamma :\mathcal{K}^{2}\to \mathcal{H}^{2}$. Below is a list of basic properties of the closure $\bar{\Gamma}$.

\begin{proposition}\label{proposition72} Let $\Gamma : \mathcal{K}^{2} \to \mathcal{H}^{2}$ be a Green's boundary relation with non-degenerate $\ran\, \Gamma$. Then:
\begin{enumerate}[(i)]%, (i), (ii),...
\item $\bar{\Gamma} :\mathcal{K}^{2} \rightarrow \mathcal{H}^{2}$ is a Green's boundary relation.
\item $ \ker \bar{\Gamma }$ is a closed symmetric relation.
\item If $\Gamma$ satisfies condition (\ref{eq222}), then also $\bar{\Gamma }$ satisfies condition (\ref{eq222}).
\item If $\Gamma$ satisfies condition (\ref{eq223}), then also $\bar{\Gamma }$ satisfies condition (\ref{eq223}).
\item If $S$ and $\hat{S}$ are closed symmetric relations associated with $\Gamma$ and $\bar{\Gamma}$, respectively, then $S=\hat{S}$. 
\item If $\hat{M}$ denotes isotropic manifold of $\dom\, \bar{\Gamma}$, then $\hat{M} \subseteq \hat{S}$. If additionally,  $\Gamma$ satisfies (\ref{eq223}), then $S=M=\hat{M} = \hat{S}$
\end{enumerate}
\end{proposition}
\textbf{Proof}. 

(i) The continuity of the scalar products (\ref{eq28}) and (\ref{eq210}) implies that the closure $\bar{\Gamma }$ satisfies (\ref{eq214}). Hence, (i) holds. 

(ii) The relation $\ker \bar{\Gamma }$ is closed as a kernel of a closed relation, and it is symmetric according to Proposition \ref{proposition26} (ii). 

(iii) The identity 
\begin{equation}
\label{eq72}
\overline{\dom}\, \bar{\Gamma }=\overline{\dom}\, \Gamma
\end{equation}
holds according to Proposition \ref{proposition34} (vi).

Now we can prove (\ref{eq222}) for $\bar{\Gamma}$. Indeed, 
\[
\left[ \hat{f},\hat{g} \right]_{\mathcal{K}^{2}}=0,\, \forall \hat{g}\in \dom\, 
\bar{\Gamma }\Rightarrow \left[ \hat{f},\hat{g} \right]_{\mathcal{K}^{2}}=0,\, \forall 
\hat{g}\in \dom\, \Gamma .
\]
According to assumption that (\ref{eq222}) holds for $\Gamma $, it follows $\hat{f}\in \overline{\dom}\, 
\Gamma =\overline{\dom}\, \bar{\Gamma }$. 

This proves (iii). 

(iv) $\left( (\ref{eq223}) \wedge (\ref{eq72})\right) \Rightarrow \left( \overline{\dom}\, \bar{\Gamma} \right)^{+}\subseteq \dom\, \Gamma \Rightarrow \left( \dom\, \bar{\Gamma} \right)^{+} \subseteq \dom\, \bar{\Gamma}$.

(v) This claim follows directly from (\ref{eq72}) and the definition of $\hat{S}=\left( \overline{\dom}\, \bar{\Gamma }\right)^{+} $. 

(vi) $\hat{M}:=\dom\, \bar{\Gamma} \cap \left( \dom\, \bar{\Gamma}\right)^{+} \subseteq \left( \dom\, \bar{\Gamma}\right)^{+} =:\hat{S}$. 

According to definitions of $M$, $\hat{M}$ and (\ref{eq72}), $M\subseteq \hat{M}$. Assume now that $\Gamma$ satisfies (\ref{eq223}). Corollary \ref{corollary56} gives $S=M$. Hence $S=M\subseteq \hat{M}\subseteq \hat{S}$. This and (v) gives $S=M= \hat{M}= \hat{S}$.
\hfill $\square$
\begin{proposition}\label{proposition73} Let $\Gamma : \mathcal{K}^{2}\to \mathcal{H}^{2}$  be a Green's boundary relation. If $\ran\,\Gamma =\mathcal{H}^{2}$, then $\Gamma $ is an operator,
\begin{equation}
\label{eq74}
\ran\, \Gamma_{i}=\Gamma_{i}\left( \ker \Gamma_{j} \right)=\mathcal{H}, \, i \neq j; \, i,j=0,1.
\end{equation}
%and
%\begin{equation}
%\label{eq76}
%\overline{\Gamma_{i}\left( \ker \Gamma_{j} \right)}=\mathcal{H}, \, i \neq j; \, i,j=0,1.
%\end{equation}
\end{proposition}
\textbf{Proof.} 
\[
\left( \Gamma^{-1}\subseteq  \Gamma^{*} \, \wedge \, \ran\,\Gamma =\mathcal{H}^{2} \right) \Rightarrow \mul\, \Gamma =  \left( \dom\, \Gamma ^{*}\right)^{[\bot]_{\mathcal{H}^{2}}} = \left( \mathcal{H}^{2}\right)^{[\bot]_{\mathcal{H}^{2}}}=\lbrace 0 \rbrace. 
\]
Hence,  $\Gamma$ is an operator. 

Let us arbitrarily select a point $h \in \ran\, \Gamma_{0}=\mathcal{H}$. Then $\left( {\begin{array}{*{20}c}
h\\
0\\
\end{array} } \right)\in \ran\, \Gamma =\mathcal{H}^{2}$. 
There exists $\hat{s}\in \dom\, \Gamma \subseteq \mathcal{K}^{2}$ such that 
$\left( {\begin{array}{*{20}c}
h\\
0\\
\end{array} } \right)=\left( {\begin{array}{*{20}c}
\Gamma_{0}\hat{s}\\
\Gamma_{1}\hat{s}\\
\end{array} } \right)$. This means that $\hat{s}\in \ker \Gamma_{1}$. i.e. $h \in \Gamma_{0}\left( \ker \Gamma_{1} \right)$. The same reasoning holds if $\Gamma_{0}$ and $\Gamma_{1}$ exchange places. Therefore (\ref{eq74}) holds. Therefore, $\ran\, \Gamma_{i}=\mathcal{H}\subseteq \Gamma_{i}\left( \ker \Gamma_{j} \right), \, i \neq j; \, i,j=0,1.$ \hfill $\square$

%If we assume $\overline{\Gamma_{0}\left( \ker \Gamma_{1} \right)}\neq \mathcal{H}$, then there exists $t_{0}\in \mathcal{H}, 0\ne t_{0}$, such that $t_{0}\left( \bot \right)_{\mathcal{H}} \overline{\ran}\, \Gamma_{0}=\overline{\Gamma_{0}\left( \ker \Gamma_{1} \right)}$. Because assumptions about $\Gamma_{0}$ and $\Gamma_{1}$ are the same, if $\ran\, \Gamma_{0}$ is not dense in $\mathcal{H}$, then also $\ran\, \Gamma_{1}$ is not dense in $\mathcal{H}$. Therefore, there exists $t_{1}\in \mathcal{H}, 0\ne t_{1}$, such that $t_{1}\left( \bot \right)_{\mathcal{H}} \overline{\ran}\, \Gamma_{1}$.
%Because $\ran\,\Gamma =\mathcal{H}^{2}$ there exists $\hat{t}$ such that $\Gamma{\hat{t}}= \left( {\begin{array}{*{20}c}
%t_{0}\\
%t_{1}\\
%\end{array} } \right)\in \mathcal{H}^{2} $. This means $t_{0} \in \ran\,\Gamma_{0}$, which contradicts $t_{0} \left( \bot \right)_{\mathcal{H}}  \overline{\ran}\, \Gamma_{0}=\overline{\Gamma_{0}\left( \ker \Gamma_{1} \right)}$. That proves (\ref{eq76}) for $i=0$ and $j=1$. When $i$ and $j$ exchange roles, we conclude that (\ref{eq76}) holds. \hfill $\square$

\begin{proposition}\label{proposition74} Let $\Gamma : \mathcal{K}^{2}\to \mathcal{H}^{2}$  be a Green's boundary relation. Assume
\begin{itemize}
\item Condition (\ref{eq222}) is fulfilled,
\item $\ran\,\Gamma =\mathcal{H}^{2}$,
\item $\dom\, \Gamma $ is closed.
\end{itemize}
Then: 
\begin{enumerate}[(i)]%, (i), (ii),...
\item Linear relations $S_{i}:=\ker \Gamma_{i},\, i=0,1$, are self-adjoint.
\item $S:=\left( \dom\, \Gamma \right)^{+}$ is a closed symmetric linear relation that satisfies $S=\ker \Gamma_{0} \cap \ker \Gamma_{1}$, i.e. extensions $S_{i}:=\ker \Gamma_{i}, i=0,1$, of $S$ are disjoint. 
\item $S^{+}=S_{0} \hat{+}S_{1}$, i.e. extensions $S_{i}, i=0,1$, of $S$ are transversal. 
\end{enumerate} 
\end{proposition}
\textbf{Proof}. According to Proposition \ref{proposition73}, $\Gamma$ is an operator. From this, and from the assumptions that $\ker \Gamma$ and $\ran\,\Gamma $ are closed, it easily follows that $\Gamma$ is a closed operator. Then, according to closed graph theorem, it follows that $\Gamma $, and hence $\Gamma_{i}, i=0,1,$ are bounded operators.

(i) The proof consists of the following four steps:
 
\textbf{Step 1.} $\ker \Gamma_{i}\subseteq \left( \ker \Gamma_{i} \right)^{+}, \, i=0,1$. This step follows from Proposition \ref{proposition26} (i).

\textbf{Step 2.} $\ran\, \Gamma =\mathcal{H}^{2}\Rightarrow \overline{\Gamma_{i}\left( \ker \Gamma_{j} \right)}=\mathcal{H},  i \neq j;\, i,\, j=0,1$.  This step follows from Proposition \ref{proposition73}.

\textbf{Step 3.} $\left( \ker \Gamma_{i} \right)^{+}\subseteq \dom\, \Gamma, \, i=0, 1$.

Indeed, because (\ref{eq222}) is assumed, we can define $S^{+}:=\overline{\dom}\, \Gamma $. According to Proposition \ref{proposition28} (i),  $S={(\overline{\dom}\, \Gamma)}^{+} \subseteq S^{+}$. 
Since $\dom\, \Gamma$ is closed and (\ref{eq222}) holds, it follows that (\ref{eq223}) holds. Additionally, $\ran\, \Gamma =\mathcal{H}^{2}$, i.e.,  $\ran\, \Gamma$ is non-degenerate. According to Corollary \ref{corollary56}, $M=N=S$, i.e., $\left(\dom\, \Gamma \right)^{+}\cap \dom\, \Gamma =\ker \Gamma =S$. 
This implies $S\subseteq \ker \Gamma_{i}$. Consequently, $ \left( \ker \Gamma_{i} \right)^{+}\subseteq S^{+}=\overline{\dom}\, \Gamma =\dom\, \Gamma $, which proves Step 3. 

\textbf{Step 4.} $\left( \ker \Gamma_{i} \right)^{+}\subseteq \ker \Gamma _{i},\, i \neq j;\, i,j=0,1$. 

Assume, for example, $\hat{f}\in \left( \ker \Gamma_{1} \right)^{+}$. According to Step 3, $\Gamma_{1}\hat{f}$ is defined. It holds 
\[
0=\left[ f',g \right]-\left[ f,g' \right]=\left( \Gamma_{1}\hat{f}, \Gamma_{0}\hat{g} \right), \forall \hat{g} \in \ker \Gamma_{1}.
\]
According to Step 2, $ \overline{\Gamma_{0}\left( \ker \Gamma_{1} \right)}=\mathcal{H}$. We conclude $\Gamma_{1}\hat{f}=0$. This proves $\left( \ker \Gamma_{1} \right)^{+}\subseteq \ker \Gamma_{1}$. By similar arguments, it holds $\left( \ker \Gamma_{0} \right)^{+}\subseteq \ker \Gamma_{0}$, which proves Step 4. 
\\

Statement (i) follows from steps 1 and 4, because  
\[
Step \,1 \wedge Step \, 4 \Rightarrow \left( \ker \Gamma_{i} \right)^{+}=\ker \Gamma_{i}, \, i=0,1.
\]

(ii) In the proof of Step 3, we saw $S=\ker \Gamma $. Because $\ker \Gamma=\ker \Gamma _{0} \cap \ker \Gamma _{1}$, statement (ii) holds. 

(iii) We will prove that every $\hat{k}=\left( {\begin{array}{*{20}c}
k\\
k'\\
\end{array} } \right)\in \dom\, \Gamma = S^{+}$ can be decomposed as $\hat{k}=\hat{u}+\hat{r}$, where $ \hat{u} \in \ker\, \Gamma _{0}$ and $\hat{r} \in \ker\, \Gamma _{1}$. 

Assume $\hat{h}:=\Gamma \hat{k}$. Then
\[
\Gamma \hat{k}=\left( {\begin{array}{*{20}c}
h\\
h'\\
\end{array} } \right)=\left( {\begin{array}{*{20}c}
0\\
h'\\
\end{array} } \right)+\left( {\begin{array}{*{20}c}
h\\
0\\
\end{array} } \right)=\Gamma \hat{t}+\Gamma \hat{r},
\]
where $ \hat{t} \in \ker\, \Gamma _{0}$ and $ \hat{r} \in \ker\, \Gamma _{1}$. This implies $ \hat{s}:=\hat{k}-\hat{t}-\hat{r} \in S \subseteq ker\, \Gamma _{0}$. Hence,
$\hat{k}:=\left( \hat{s}+ \hat{t}\right)+\hat{r} =: \hat{u}+ \hat{r} \in \ker\, \Gamma _{0} \hat{+} \ker\, \Gamma _{1}$. This proves (iii). \hfill $\square$
\\

According to Corollary \ref{corollary58}, $\Gamma$ is an ordinary boundary relation if and only if the conditions of Proposition \ref{proposition74} are satisfied. This implies that claim (iii) of Proposition \ref{proposition74} could be established through Corollary \ref{corollary58} and \cite[Proposition 3.4]{BDHS}. However, the proof of Proposition \ref{proposition74} (iii) that we presented is elementary and direct, not based on the corresponding statement in Hilbert spaces and their fundamental symmetries. It also demonstrates how Green's boundary relation approach can be employed.
\section{Fitting of generalized boundary triples into Green's boundary model in the Krein space}\label{s16}

When $\mathcal{K}$ is a separable Hilbert space, and $S\subseteq \mathcal{K}^{2}$ is a closed symmetric linear operator with equal finite or infinite deficiency indexes, various kinds of boundary triples $\Pi =\left(\mathcal{H},\Gamma_{0}, \Gamma_{1} \right)  $, with  $\mathcal{K}$ being a Hilbert space and $\Gamma : \mathcal{K}^{2}\to \mathcal{H}^{2}$ being an operator, were defined in \cite{BL,DHM}: ordinary, B-generalized, AB-generalized, isometric, unitary, S-generalized, quasi boundary triple, and others. 

In this section, we will fit all the aforementioned boundary triples into the Green's boundary relation model and generalize them to a Krein space $\mathcal{K}$ and linear relation $S$ and $\Gamma$ using the Green's boundary relation model.

The definitions of isometric and unitary boundary relations in terms of Green's boundary relations were provided earlier; see definitions \ref{definition210} and \ref{definition214}. These definitions also extend the concepts of isometric and unitary triples given in \cite[Definition 1.9]{DHM} to Krein spaces and linear relations $S$ and $\Gamma $.

The definitions of AB-generalized and B-generalized boundary triples, as seen in \cite[definitions 1.5 and 1.8]{DHM} and papers referenced therein, can be generalized in terms of Green's boundary relation in the Krein space $\mathcal{K}$ as follows:

\begin{definition}\label{definition82} A Green's boundary relation $\Gamma : \mathcal{K}^{2}\to \mathcal{H}^{2}$ is called an almost B-generalized, or AB-generalized boundary relation if it satisfies:

\begin{itemize}
\item Condition (\ref{eq222}),
\item $\overline{\ran}\, \Gamma_{0}=\mathcal{H}$, and
\item $A:=\ker \Gamma_{0}$ is a self-adjoint linear \textbf{relation} in $\mathcal{K}$.
\end{itemize}
If additionally ${\ran\, \Gamma }_{0}=\mathcal{H}$, then the relation $\Gamma$ is called a B-generalized boundary relation.
\end{definition}

It is easy to see that the following inclusions hold: 
\\

$Ordinary \, Boundary \, Relations \subseteq B$-$boundary  \, Relations \subseteq AB$-$boundary  \, Relations $
\[
 \subseteq  Isometric \, Boundary \, Relations \subseteq Green's \, Boundary \, Relations.
\]

For the definition of a \textit{quasi boundary triple} for $S^{+}$ in the Hilbert space $\mathcal{K}$, refer to \cite[Definition 2.1]{BL}. Here is a generalization of the definition to the Krein space $\mathcal{K}$ in terms of Green's boundary relation.

\begin{definition}\label{definition810} A Green's boundary relation $\Gamma : \mathcal{K}^{2}\to \mathcal{H}^{2}$ is called quasi-boundary if it satisfies: 

\begin{itemize}
\item Condition (\ref{eq222}),
\item $\overline{\ran}\, \Gamma =\mathcal{H}^{2},$ \textit{and }
\item $A\, :=\ker \Gamma_{0}$ is a self-adjoint operator in $\mathcal{K}$.
\end{itemize}
\end{definition}
Note that a quasi-boundary relation $\Gamma$ is always a closable isometric boundary operator.

\begin{definition}\label{definition812} A Green's boundary relation $\Gamma : \mathcal{K}^{2}\to \mathcal{H}^{2}$ is called an S-generalized boundary relation if:

\begin{itemize}
\item It is a unitary relation,
\item $A:=\ker \Gamma_{0}$ is a self-adjoint operator.
\end{itemize}
\end{definition}

\textbf{Compliance with Ethical Standards Statements:}

\textbf{Conflict of Interest:} The author declares that there is no conflict of interest. 

\textbf{Funding:} No funding was received to assist with the preparation of this manuscript.

\textbf{Ethical Conduct:}  This research is original. The manuscript is being submitted only to Journal of Mathematical Sciences. All previous results used in this paper are appropriately cited. 

\textbf{Data Availability Statements:} This research is purely mathematical. The data supporting this paper consist solely of mathematical statements, which are openly available in the cited references. There are no other data involved besides mathematical statements.

%BML, Actuarial Department

%120 Royall Street. 

%Canton, MA 02131,USA

%muhamed.borogovac@gmail.com 

\end{document}